\documentclass[a4paper,12pt]{article}
\usepackage{comment,amsmath,amssymb,amsthm,pstricks,changepage,pifont,hyperref,pstricks-add}
\usepackage[english]{babel}

\newtheorem{theorem}{Theorem}[section]

\newtheorem{lemma}[theorem]{Lemma}
\newtheorem{corollary}[theorem]{Corollary}

\newtheorem{conjecture}[theorem]{Conjecture}
\theoremstyle{remark}

\DeclareMathOperator{\lsize}{ls}
\DeclareMathOperator{\lwidth}{lw}
\DeclareMathOperator{\conv}{conv}
\DeclareMathOperator{\Vol}{Vol}

\begin{document}

\title{The lattice size of a lattice polygon}
\author{Wouter Castryck and Filip Cools}

\maketitle

\begin{abstract}
  \noindent We give upper bounds on the minimal degree of a model in $\mathbb{P}^2$ and the minimal bidegree of a model
  in $\mathbb{P}^1 \times \mathbb{P}^1$
  of the curve defined by a given Laurent polynomial, in terms of the combinatorics of the Newton polygon of the latter. 
  We prove in various cases that this bound is 
  sharp as soon as the polynomial is sufficiently generic with respect to its Newton polygon.\\
   
\noindent \emph{MSC2010:} Primary 14H45, Secondary 14H51, 14M25
\end{abstract}

\section{Introduction}

Let $k$ be an algebraically closed field and
let $f \in k[x^{\pm 1}, y^{\pm 1}]$ be an irreducible Laurent polynomial,
whose Newton polygon we denote by $\Delta(f)$ (which we assume to be two-dimensional).
Let $\mathbb{T}^2 = k^* \times k^*$ be the two-dimensional torus over $k$, and denote by $U_f \subset \mathbb{T}^2$ the
curve defined by $f$. (Throughout this paper, all curves
are understood to be irreducible, but not necessarily non-singular
and/or projective.)
For a curve $C/k$ we define $s_2(C)$ as the minimum of 
    \begin{align*}
      S_2(C) & = \left\{ \, \left. d \in \mathbb{N} \, \right| \, \text{$C \simeq$ a curve of degree $d$ in $\mathbb{P}^2$} \, \right\}  \\
      \intertext{and $s_{1,1}(C)$ as the lexicographic minimum of} 
      S_{1,1}(C) & = \left\{ \, \left. (a,b) \in \mathbb{N}^2 \, \right| \, \text{$a \leq b$ and $C \simeq$ a curve of bidegree $(a,b)$ in $\mathbb{P}^1 \times \mathbb{P}^1$} \, \right\},
    \end{align*}
where $\simeq$ denotes birational equivalence. The aim of this article is to give upper bounds on the invariants
$s_2(U_f)$ and $s_{1,1}(U_f)$ purely in terms of the combinatorics of $\Delta(f)$.

The invariant $s_2(C)$ has seen study in the past \cite{CoppensMartens,harui,KeemMartens} but is not well-understood.
On the other hand we are unaware of existing literature explicitly devoted to $s_{1,1}(C)$, even though 
for hyperelliptic curves the notion did make an appearance \cite{freykani} in the context of 
cryptography. Note that at first sight, the definition of $s_{1,1}(C)$ has a non-canonical flavor: instead of lexicographic, one could also consider the minimum with respect to
other types of monomial orders on $\mathbb{N}^2$. But in fact we conjecture: 
\begin{conjecture} \label{conjecturebidegree} For each curve $C/k$ the set $S_{1,1}(C)$
admits a minimum with respect to the product order $\leq \times \leq$ on $\mathbb{N}^2$. 
\end{conjecture}
\noindent Because the product order is coarser than every monomial order, this would mean
that the term `lexicographic' can be removed without ambiguity. In Section~\ref{section_basicfacts}
we will state a number of basic facts on $s_2(C)$ and $s_{1,1}(C)$, along with some motivation
in favor of Conjecture~\ref{conjecturebidegree}. 

Our central combinatorial notion is the \emph{lattice size} $\lsize_X(\Delta)$ of a lattice polygon $\Delta$ with respect to
a set $X \subset \mathbb{R}^2$ with positive Jordan measure. 
In case $\Delta \neq \emptyset$ we
define it as the smallest integer $d \geq 0$ for which there exists a unimodular 
transformation $\varphi : \mathbb{R}^2 \rightarrow \mathbb{R}^2$ such that
\[ \varphi(\Delta) \subset d X.\]
A unimodular transformation that attains this minimum is said to compute the lattice size.
We will restrict ourselves to three instances of $X$, namely
\[ [0,1] \times \mathbb{R}, \quad \ \Sigma = \conv \{(0,0), (1,0), (0,1) \}, \quad \ \square =  \conv \{ (0,0), (1,0), (0,1), (1,1) \}, \]
where it is convenient to define $\lsize_X(\emptyset)= -1, -2, -1$, respectively.

In the case of $X = \Sigma$ the lattice size measures the smallest standard triangle containing a unimodular copy of $\Delta$.
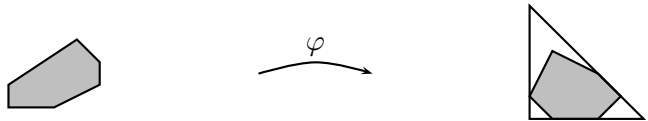
\begin{center}
\psset{unit=0.3cm}
\begin{pspicture}(0,-0.5)(5,5)
\pspolygon[fillstyle=solid,fillcolor=lightgray](0,0)(0,1)(3,3)(4,2)(4,1)(2,0)
\end{pspicture}
\qquad \qquad
\begin{pspicture}(0,0)(5,5)
\pscurve{->}(0,2)(2.5,2.5)(5,2)
\rput(2.5,3.2){\small $\varphi$}
\end{pspicture}
\qquad \qquad \ \
\psset{unit=0.3cm}
\begin{pspicture}(0,0)(5,5)
\pspolygon[fillstyle=solid,fillcolor=lightgray](1,0)(0,1)(1,3)(3,2)(4,1)(3,0)
\pspolygon(0,0)(5,0)(0,5)
\end{pspicture}
\end{center}
This was studied by Schicho \cite{schicho}, who designed an algorithm for
finding a unimodular transformation that maps a given polygon $\Delta$ inside a small standard triangle. He did this in the context of simplifying parametrizations of rational surfaces. Our results below show
that Schicho's algorithm works optimally, that is, its output computes the lattice size $\lsize_\Sigma(\Delta)$. 
In the case of $X = [0,1] \times \mathbb{R}$ the lattice size is nothing else than the commonly studied lattice width, which we
denote by $\lwidth(\Delta)$ rather than $\lsize_{[0,1] \times \mathbb{R}}(\Delta)$. See \cite[Lem.\,5.2]{linearpencils} for some of its properties, such 
as Fejes T\'oth and Makai Jr.'s result that $\lwidth(\Delta)^2 \leq 8 \text{Vol}(\Delta) / 3$.
In the case of $X = \square$ the notion implicitly appears in the work of Arnold \cite{arnold} and Lagarias--Ziegler \cite[Thm.\,2]{lagariasziegler} in the context of counting
lattice polygons (up to unimodular equivalence) with a given volume; they found that $\lsize_\square(\Delta) \leq 4 \text{Vol}(\Delta)$ as soon as $\Delta$ is two-dimensional. Note that this implies
the bound $\lsize_\Sigma(\Delta) \leq 8 \text{Vol}(\Delta)$, which is most likely not sharp.

Recently, Lubbes and Schicho \cite[Thm.\,13]{lubbesschicho} and the current authors~\cite[Thm.\,4]{curvegonalities} independently provided an explicit formula for $\lwidth(\Delta)$
in terms of $\lwidth(\Delta^{(1)})$, where $\Delta^{(1)}$ denotes the convex hull of the lattice points
in the interior of $\Delta$; see Lemma~\ref{lubbesschicholemma} for a precise statement. This yields 
a recursive method for computing the lattice width in practical situations,
by gradually `peeling off' the polygon.\footnote{We remark that
for very large polygons there exist more effective methods; see e.g.\ \cite{freschet}.}
The biggest part of this article (Sections~\ref{section_recursivesigma} and~\ref{section_recursivesquare}) is devoted to proving similar recursive formulas for $\lsize_\Sigma(\Delta)$ and $\lsize_\square(\Delta)$, which
can again be used for computing the lattice size in practice. 
In the former case one recovers Schicho's algorithm.
In the latter case the proof entails that the unimodular transformations computing $\lsize_\square(\Delta)$ essentially also compute $\lwidth(\Delta)$.
This is made precise in Section~\ref{section_minproductorder}, where as a corollary we obtain:
\begin{theorem} \label{minproductorder}
For each non-empty lattice polygon $\Delta$ the set
\begin{align*} 
    S_{1,1}(\Delta) & = \left\{ \, \left. (a,b) \in \mathbb{N}^2 \, \right| \, \text{$a \leq b$ and $\exists \, \Delta' : \Delta \simeq \Delta'$ with $\Delta' \subset [0,a] \times [0,b]$} \, \right\}
\end{align*}
admits a minimum with respect to the product order on $\mathbb{N}^2$, namely $s_{1,1}(\Delta) := (\lwidth(\Delta),\lsize_\square(\Delta))$.
\end{theorem}
\noindent Here $\simeq$ denotes unimodular equivalence. We will sometimes (but not always) write
$\square_{a,b}$ instead of $[0,a] \times [0,b]$. The reader can view
Theorem~\ref{minproductorder} as a combinatorial version of Conjecture~\ref{conjecturebidegree}.

Now if we write
\[ f = \sum_{(i,j) \in \mathbb{Z}^2} c_{i,j} x^i y^j \quad \in k[x^{\pm 1}, y^{\pm 1}] \]
then for every unimodular transformation $\varphi : \mathbb{R}^2 \rightarrow \mathbb{R}^2 $ 
the Laurent polynomial
\[ f^\varphi = \sum_{(i,j) \in \mathbb{Z}^2} c_{i,j} x^{\varphi_1(i,j)} y^{\varphi_2(i,j)} \]
(with $\varphi_1$ and $\varphi_2$ the component functions of $\varphi$) satisfies $\Delta(f^\varphi) = \varphi(\Delta(f))$. Since $U_f$ and $U_{f^\varphi}$ are isomorphic
 it follows that 
\begin{equation} \label{trivialbounds}
 s_2(U_f) \leq \lsize_\Sigma(\Delta(f)) \quad \text{and} \quad s_{1,1}(U_f) \leq s_{1,1}(\Delta(f)),
\end{equation}
where the second inequality should be read lexicographically. While the first bound is straightforward, we note that the second bound 
relies on Theorem~\ref{minproductorder}.
Our main result, which shows up as a consequence to our recursive formulas for the lattice size, 
refines these bounds:

\begin{theorem} \label{s2upper}
 One has
 \begin{equation} \label{tobeprovenupperbounds}
   s_2(U_f) \leq \lsize_\Sigma(\Delta(f)^{(1)}) + 3 \quad \text{and} \quad s_{1,1}(U_f) \leq s_{1,1}(\Delta(f)^{(1)}) + (2,2). 
 \end{equation}
 If $\Delta(f) \simeq d\Upsilon$ for some $d \geq 2$ then the first bound sharpens to $s_2(U_f) \leq 3d-1$;
 if $d=2$ then also the second bound sharpens to
 $s_{1,1}(U_f) \leq (3,4)$.
\end{theorem}

\noindent Here $\Upsilon = \conv \{(-1,-1),(1,0),(0,1) \}$. The proof of Theorem~\ref{s2upper} is given in Section~\ref{section_maintheorem}. We will see below that $\lsize_\Sigma(\Delta^{(1)}) + 3 \leq \lsize(\Delta)$ and $s_{1,1}(\Delta^{(1)}) + (2,2) \leq s_{1,1}(\Delta)$ 
as soon as $\Delta$ is two-dimensional,
and that the difference can be arbitrarily large.
Thus Theorem~\ref{s2upper} can be seen as a considerable improvement over the bounds (\ref{trivialbounds}).
As a teasing example, consider a hyperelliptic curve $C$ of genus $g \geq 2$ defined by a Weierstrass equation
\[ f := y^2 + h_1(x)y + h_2(x) = 0, \]
with $h_2 \in k[x]$ of degree $2g + 1$ and $h_1 \in k[x]$ of degree at most $g$. Assume for simplicity that $h_2(0) \neq 0$, so
that the Newton polygon $\Delta(f)$ equals
  \begin{center}
  \psset{unit=0.3cm}
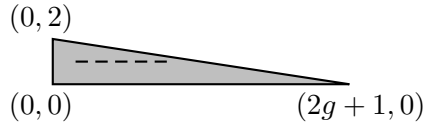
  \begin{pspicture}(-1,-1.3)(14,3.3)
  \pspolygon[fillstyle=solid,fillcolor=lightgray](0,0)(13,0)(0,2)
  \psline[linestyle=dashed](1,1)(5,1)
  \rput(-0.5,-1){\small $(0,0)$}
  \rput(13.5,-1){\small $(2g+1,0)$}  
  \rput(-0.5,2.9){\small $(0,2)$}  
  \end{pspicture}
  \end{center}
The interior polygon $\Delta(f)^{(1)}$ equals $\conv \{ (1,1), (g,1) \}$; it is indicated by the dashed line.
In this case the bounds (\ref{trivialbounds}) read $s_2(C) \leq 2g+1$ and $s_{1,1}(C) \leq (2,2g+1)$, while
Theorem~\ref{s2upper} yields $s_2(C) \leq g+2$ and $s_{1,1}(C) \leq (2,g+1)$. The latter bounds are actually sharp; see Section~\ref{section_basicfacts}.
More generally, we conjecture: 
\begin{conjecture} \label{sufficientlygenericconjecture}
 If $f$ is sufficiently generic with respect to its Newton polygon $\Delta(f) \not \simeq 2\Upsilon$, then 
the (smallest applicable) bounds of Theorem~\ref{s2upper} are met. 
\end{conjecture}
\noindent In Section~\ref{section_conjecture},
where we will be more precise on what is meant by `sufficiently generic', we will prove
this conjecture in a number of special cases.

\subsection*{Acknowledgments} We sincerely thank the anonymous referees, whose comments on a first version of this article
led to many extensions and improvements. We are also grateful to Marc Coppens for answering several of our 
questions on linear systems. 
This research was conducted in the framework of Research Project G093913N of the Research Foundation - Flanders (FWO).

\section{Basic facts on the minimal (bi)degree} \label{section_basicfacts}

Let $C$ be a curve of (geometric) genus $g$ over an algebraically closed field $k$.
In this section we discuss a number of basic properties of the invariants $s_2(C)$ and $s_{1,1}(C)$.
Throughout we make the assumption that $\text{char} \, k = 0$, because several of our references rely on it.

In the case of $s_2(C)$ it is known that
\[ \frac{3 + \sqrt{8g + 1}}{2} \leq s_2(C) \leq g + 2.\]
The lower bound is met if and only if $C$ is birationally
equivalent to a non-singular projective plane curve. As for the upper bound one has $s_2(C) = g + 2$ if and only if
$C$ is elliptic or hyperelliptic. If $g \geq 6$ then $s_2(C) = g + 1$ if and only if $C$ is bi-elliptic.
See \cite{harui} and the references therein for proofs.

In the case of $s_{1,1}(C)$ we prove an analogous statement:
\begin{lemma} \label{marclemma}
One has $s_{1,1}(C) = (c,d)$, where $c$ is the gonality of $C$ and $d$ satisfies
\[ \frac{g}{c-1} + 1 \leq d \leq g + 1 \]
unless $c = 1$ where $d = 1$. The lower bound is met if and only if $C$ is birationally
equivalent to a non-singular curve in $\mathbb{P}^1 \times \mathbb{P}^1$. If $C$ is elliptic or hyperelliptic then the upper bound is met.
\end{lemma}

\noindent \textsc{Proof.} If $c=1$ then the statement is trivial, so we can assume that $c \geq 2$.

As for the upper bound, fix a $g^1_c$ on $C$ and pick
a point $P \in C$. Let $D \in g^1_c$ be such that $P$ is in the support.
Now construct a divisor $D'$ by gradually adding points that are not in the support of $D$ to the point $P$, until $\dim H^0(C,D') = 2$.
By the Riemann-Roch theorem this happens after at most $g$ steps, i.e.\ $d := \deg D' \leq g + 1$. By construction, 
the corresponding base-point free $g^1_d$ does not have a factor in common with our given $g^1_c$, so
we can use $g^1_c \times g^1_d$ to map $C$ to a birationally equivalent curve of bidegree $(c,d)$.

As for the other inequality, consider Baker's bound~\cite{beelen}, which 
says that the genus of the curve defined by an irreducible Laurent polynomial $f \in k[x^{\pm 1},y^{\pm 1}]$ is bounded by
$\sharp (\Delta(f)^{(1)} \cap \mathbb{Z}^2)$. Now the Newton polygon of a polynomial of bidegree $(c,d)$ is contained
in the rectangle:
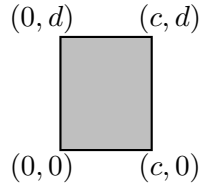
\begin{center}
\psset{unit=0.3cm}
\begin{pspicture}(-1,-1)(5,6)
\pspolygon[fillstyle=solid,fillcolor=lightgray](0,0)(4,0)(4,5)(0,5)
\rput(-0.8,-0.8){\small $(0,0)$}
\rput(4.8,-0.8){\small $(c,0)$}
\rput(4.8,5.8){\small $(c,d)$}
\rput(-0.8,5.8){\small $(0,d)$}
\end{pspicture}
\end{center}
Hence $g \leq (c-1)(d-1)$, from which the lower bound follows. If
there is a singularity in $\mathbb{P}^1 \times \mathbb{P}^1$, then without loss of generality we may assume that it concerns an affine point $(x_0,y_0)$. But 
then the Newton polygon of $f(x + x_0, y + y_0)$ is contained in:
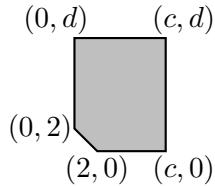
\begin{center}
\psset{unit=0.3cm}
\begin{pspicture}(-1,-1)(5,6)
\pspolygon[fillstyle=solid,fillcolor=lightgray](0,1)(1,0)(4,0)(4,5)(0,5)
\rput(1,-0.8){\small $(2,0)$}
\rput(-1.5,1){\small $(0,2)$}
\rput(4.8,-0.8){\small $(c,0)$}
\rput(4.8,5.8){\small $(c,d)$}
\rput(-0.8,5.8){\small $(0,d)$}
\end{pspicture}
\end{center}
Therefore $g \leq (c-1)(d-1)-1$, which shows that the lower bound cannot be attained in this case.

Finally, if $C$ is elliptic or hyperelliptic then $c = 2$ so that the lower bound meets the upper bound. 
\hfill $\blacksquare$\\

We think that for $c \geq 2$ the upper bound can be improved to $g + 3 - c$. 
Namely, by Brill-Noether theory the space of $g^1_{g + 3 - c}$'s on $C$ has dimension $g - 2c + 4$, while the subspace
of pencils of the form $g^1_c + \text{base points}$ has dimension $g - 2c + 3$. This gives plenty of base-point free
$g^1_d$'s with $d \leq g + 3 - c$ that do not obviously have a factor in common with the given $g^1_c$. But 
we did not succeed in proving that there indeed always exists such a truly independent $g^1_d$. 
The bi-elliptic case illustrates the subtlety of the argument: here one has a full-dimensional component of dependent $g^1_{g + 3 - c}$'s. Nevertheless the bound $g + 3 - c = g - 1$ is valid here (and met); see \cite[Ex.\,1.13]{CoppensMartens}.

As a special cases of Conjecture~\ref{conjecturebidegree}, we note:

\begin{lemma} \label{conjecturebidegreespecialcases}
If the gonality of $C$ is a prime number 
then $S_{1,1}(C)$ admits a minimum with respect to the product order $\leq \times \leq$ on $\mathbb{N}^2$.
\end{lemma}

\noindent \textsc{Proof.} Fix a gonality pencil $g^1_c$. It suffices to show that if $(a,b) \in S_{1,1}(C)$ then $(c,a) \in S_{1,1}(C)$ or $(c,b) \in S_{1,1}(C)$.
In other words, it is sufficient to prove that at least one of the given $g^1_a$ and $g^1_b$ is independent of our $g^1_c$. But if $g^1_a$ and $g^1_c$ have
a common factor, then by primality this factor must be $g^1_c$ itself, and similarly for $g^1_b$. Because $g^1_a$ and $g^1_b$ are mutually independent, the claim follows.
\hfill $\blacksquare$\\

We do not have much additional evidence in favor of Conjecture~\ref{conjecturebidegree}, except that 
all our attempts to construct a counterexample failed in a suspicious way: each time unexpected linear pencils
popped up that made the statement true. As a typical example, we considered the fiber product  
\[ C : \left\{ \begin{array}{l}  y_1^3 = f_1(x) \\  y_2^3 = f_2(x) \\ \end{array} \right. \]
of two cyclic degree $3$ covers
of the projective line, with $f_1(x), f_2(x)$ degree $6$ polynomials that are squarefree and mutually coprime.
This is a $9$-gonal curve of genus $28$ by Riemann-Hurwitz, so in view of Lemma~\ref{marclemma} we have $(9,d) \in S_{1,1}(C)$ with $d \leq 29$. 
On the other hand both covers naturally admit a $g^1_4$, which when composed with the $g^1_3$ of
the other curve gives rise to
two independent $g^1_{12}$'s on $C$, each of which has a component in common
with our $g^1_9$. So we also find that $(12,12) \in S_{1,1}(C)$, while it
is not obvious that $(9,e) \in S_{1,1}(C)$ with $e \leq 12$, especially because the genus is so high. However, in all
concrete versions that we tried the substitution $y_1 \leftarrow y_1 + y_2$, when followed by a projection on the $(x,y_1)$-plane,
resulted in a plane degree $15$ curve having several triple points, each of which corresponds to a $g^1_{12}$ by projection. In this way we
always found that $(9,12) \in S_{1,1}(C)$.

\section{A recursive formula for $\lsize_\Sigma(\Delta)$} \label{section_recursivesigma}

We begin by discussing some first properties. For $d \in \mathbb{Z}_{\geq 0}$ one has
\[ \lsize_\Sigma \left( \conv \{(0,0), (d,0)\} \right) = d. \]
Indeed, it is immediate that 
$\conv \{(0,0), (d,0)\} \subset d\Sigma$
and
that the integral distance $\gcd (a_2-a_1,b_2-b_1)$ between two points $(a_1,b_1)$, $(a_2,b_2) \in (d-1)\Sigma$
cannot exceed $d-1$. More generally, every lattice polygon that contains a line segment of integral length $d$
must have lattice size at least $d$ with respect to $\Sigma$.
In particular $\lsize_\Sigma(d\Sigma) = d$.

\begin{lemma} \label{spectrum}
  Let $\Delta$ be a non-empty lattice polygon. Then $\lwidth(\Delta) \leq \lsize_\Sigma(\Delta)$, and equality holds if and only if $\Delta \simeq d\Sigma$
  for some integer $d \geq 0$.
\end{lemma}

\noindent \textsc{Proof.} This follows because $\lwidth(d\Sigma) = d$, while every strict subpolygon $\Gamma \subset d\Sigma$
satisfies $\lwidth(\Gamma) < d$. 
\hfill $\blacksquare$\\

A less straightforward lattice size calculation is:

\begin{lemma} \label{latticesizesquare}
Let $a,b \in \mathbb{Z}_{\geq 0}$ and consider 
$\square_{a,b} = [0,a] \times [0,b]$.
 Then $\lsize_\Sigma(\square_{a,b}) = a + b$.
\end{lemma}
\noindent \textsc{Proof.} 
The case where $a = 0$ or $b = 0$ follows from the above considerations, so we can 
assume that $a, b \geq 1$. Instead of looking for the minimal $d$ such that $\square_{a,b}$ can be mapped 
inside $d\Sigma$ through a unimodular transformation, we will look for the minimal $d$ such that $\square_{a,b}$
is contained in a unimodular transform of $d\Sigma$. More precisely, we will prove the following assertion by induction on $a+b$: 
\begin{itemize}
 \item[] \emph{We have
$\lsize_\Sigma(\square_{a,b}) = a + b$. Moreover, there are exactly four ways of fitting
$\square_{a,b}$ inside a unimodular transform of $(a+b)\Sigma$:}
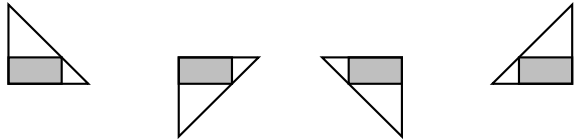
\begin{center}
\psset{unit=0.35cm}
\begin{pspicture}(-2,-3)(4,3)
\pspolygon[fillstyle=solid,fillcolor=lightgray](0,0)(2,0)(2,1)(0,1)
\pspolygon(0,0)(3,0)(0,3)
\end{pspicture}
\begin{pspicture}(-2,-3)(4,3)
\pspolygon[fillstyle=solid,fillcolor=lightgray](0,0)(2,0)(2,1)(0,1)
\pspolygon(0,1)(3,1)(0,-2)
\end{pspicture}
\begin{pspicture}(-2,-3)(4,3)
\pspolygon[fillstyle=solid,fillcolor=lightgray](0,0)(2,0)(2,1)(0,1)
\pspolygon(2,1)(2,-2)(-1,1)
\end{pspicture}
\begin{pspicture}(-2,-3)(4,3)
\pspolygon[fillstyle=solid,fillcolor=lightgray](0,0)(2,0)(2,1)(0,1)
\pspolygon(2,0)(2,3)(-1,0)
\end{pspicture}
\end{center}
\end{itemize}
The basis of our induction is the case $a=b=1$. Here, the first part of the assertion holds
because
$\square_{1,1} \subset 2\Sigma$ and $\Vol(\square_{1,1}) > \Vol(\Sigma)$. 
The second part follows because $2\Sigma$ contains only $3$ lattice points that are non-vertices. Therefore, when
fitting $\square_{1,1}$ inside a transform of $2\Sigma$, at least
one of its vertices must coincide with a vertex of $2\Sigma$, and the two adjacent
vertices of $\square_{1,1}$ must coincide with the interior lattice points of the respective adjacent edges of $2\Sigma$. 
From this the claim follows easily. 

Now assume that $a, b \geq 1$ and (without loss of generality) that $a \geq 2$. Clearly $\square_{a,b} \subset (a + b)\Sigma$. 
Suppose that $\square_{a,b}$ sits inside a unimodular transform
of $(a + b - 1)\Sigma$. By applying the induction hypothesis to 
$\square_{a-1,b} \subset \square_{a,b}$ 
we find that $(a + b - 1)\Sigma$ must enclose this subpolygon in one of the four manners
above. But for each of these four configurations, it is clear that $\square_{a,b}$ itself could not have
been contained in $(a + b - 1)\Sigma$: contradiction.
As for the second assertion, let $\Sigma'$ be a unimodular transform of $(a + b)\Sigma$ containing
$\square_{a,b}$. Then
\begin{itemize}
 \item each edge of $\Sigma'$ must contain at least one vertex of $\square_{a,b}$: otherwise
we could crop $\Sigma'$ to a unimodular transform of $(a+b-1)\Sigma$ that still contains $\square_{a,b}$;
 \item at least one vertex $v$ of $\Sigma'$ does not appear as a vertex of $\square_{a,b}$: otherwise
the latter would be a triangle; 
 \item the edges of $\Sigma'$ that are adjacent to $v$ cannot contain two vertices of $\square_{a,b}$ each: otherwise
$\square_{a,b}$ would contain two non-adjacent non-parallel edges.
\end{itemize}
So there must be an edge $\tau \subset \Sigma'$ that contains exactly one vertex $v$ of $\square_{a,b}$. 
Then the transform of $(a + b - 1)\Sigma$
obtained from $\Sigma'$ by shifting $\tau$ inwards contains $\square_{a,b} \setminus \{v\}$. In particular
it contains (a translate of) $\square_{a-1,b}$. By applying the induction hypothesis we
find that $\Sigma'$ must be positioned in one of the four standard ways above.
\hfill $\blacksquare$\\

 We now investigate the relation between $\lsize_\Sigma(\Delta)$ and $\lsize_\Sigma(\Delta^{(1)})$.
Since for $d \geq 3$ one has $(d\Sigma)^{(1)} \simeq (d-3)\Sigma$, we have that
\begin{equation} \label{latticesize_trivialbound}
 \lsize_\Sigma(\Delta^{(1)}) \leq \lsize_\Sigma(\Delta) - 3
\end{equation}
 as soon as $\Delta$ is two-dimensional (this includes the case where $\Delta^{(1)} = \emptyset$, which can be verified
separately).
Typically, one expects equality to hold, but there are many exceptions, which are classified by Theorem~\ref{latticesize_classification} below.
%

In what follows, we will make use of the following terminology and facts;
see \cite[\S4]{HaaseSchicho} or \cite[\S2.2]{Koelman} for proofs.
An edge $\tau$ of a two-dimensional lattice polygon $\Gamma$ is
always supported on a line $a_\tau X + b_\tau Y = c_\tau$ with $a_\tau, b_\tau, c_\tau \in \mathbb{Z}$ and
$a_\tau, b_\tau$ coprime. When signs are chosen appropriately, we can
moreover assume that $\Gamma$ is contained in the half-plane $a_\tau X + b_\tau Y \leq c_\tau$. 
The line $a_\tau X + b_\tau Y = c_\tau + 1$ is called the \emph{outward shift} of $\tau$. It is denoted
by $\tau^{(-1)}$, and the polygon (which may take vertices outside $\mathbb{Z}^2$) 
that arises as the intersection of the half-planes $a_\tau X + b_\tau Y \leq c_\tau + 1$
is denoted by $\Gamma^{(-1)}$. 
If $\Gamma = \Delta^{(1)}$ for some lattice polygon $\Delta$, then the outward shifts
of two adjacent edges of $\Gamma$ always intersect in a lattice point, and in fact $\Gamma^{(-1)} = \Delta^{(1)(-1)}$
is a lattice polygon. Moreover, $\Delta \subset \Delta^{(1)(-1)}$, i.e.\ $\Delta^{(1)(-1)}$ is the maximal 
lattice polygon (with respect
to inclusion) for which the convex hull of the interior lattice points equals $\Delta^{(1)}$.

Before stating Theorem~\ref{latticesize_classification}, let us prove two auxiliary lemmas:

\begin{lemma} \label{latticesize_auxiliary1}
Assume that there exist parallel edges
$\tau \subset \Delta$ and $\tau' \subset \Delta^{(1)}$ whose
supporting lines are at integral distance $1$ of each other, 
of respective lengths $r$ and $s$. If $r \geq s + 3$ then
$\lsize_\Sigma(\Delta^{(1)}) = s$ and $\lsize_\Sigma(\Delta) = r$.
\end{lemma}

\noindent \emph{Remark.} As usual, by an edge we mean a one-dimensional face. 
In particular, if $\Delta^{(1)}$ is one-dimensional
then it is an edge of itself. Example: consider the hyperelliptic Weierstrass polygon
\[ \conv \{ (0,0), (2g+1,0), (0,2) \} \]
from the introduction.
Then $\lsize_\Sigma(\Delta^{(1)}) = g$ and $\lsize_\Sigma(\Delta) = 2g + 1$. This shows that the difference 
between $\lsize_\Sigma(\Delta)$ and $\lsize_\Sigma(\Delta^{(1)})$ can be
arbitrarily large.\\

\noindent \textsc{Proof of Lemma~\ref{latticesize_auxiliary1}.} 
By using a unimodular transformation if needed, we can assume that $\tau = \conv \{(-1,-1), (r-1,-1)\}$
and $\tau' = \conv \{(0,0), (s,0)\}$. Since $r \geq s + 3$ and $\Delta^{(1)}$ cannot
contain any lattice points on the line $Y = 0$ apart from those contained in $\tau'$,
\begin{itemize}
 \item the edge of $\Delta$ that is left-adjacent to $\tau$ must pass through or to the right of $(-1,0)$, and
 \item the edge of $\Delta$ that is right-adjacent to $\tau$ must pass through or to the left of $(s+1,0)$.
\end{itemize}
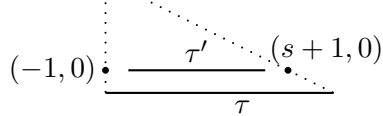
\begin{center}
\psset{unit=0.3cm}
\begin{pspicture}(-2, -1.6)(9,3)
\psline(-1,-1)(9,-1)
\psline(0,0)(6,0)
\rput(5,-1.6){\small $\tau$}
\rput(3,0.85){\small $\tau'$}
\pscircle[fillstyle=solid,fillcolor=black](-1,0){0.15}
\pscircle[fillstyle=solid,fillcolor=black](7,0){0.15}
\psline[linestyle=dotted](-1,-1)(-1,3)
\psline[linestyle=dotted](9,-1)(1,3)
\rput(-3.3,0){\small $(-1,0)$}
\rput(8.7,0.9){\small $(s + 1,0)$}
\end{pspicture}
\end{center}
From the convexity of $\Delta$ one immediately sees that
$\Delta \subset (-1,-1) + r\Sigma$, and similarly that $\Delta^{(1)} \subset s \Sigma$.
Therefore $\lsize_\Sigma(\Delta^{(1)}) \leq s$ and $\lsize_\Sigma(\Delta) \leq r$, and equality follows from
the considerations preceding Lemma~\ref{latticesizesquare}.
\hfill $\blacksquare$

\begin{lemma} \label{latticesize_auxiliary2}
Assume that $\Delta^{(1)}$ is two-dimensional.
Let $s \geq 1$ be an integer such that $\Delta^{(1)} \subset s\Sigma$, and assume 
that $\Delta^{(1)}$ has an edge $\tau'$ 
in common with $s\Sigma$. Let $\tau'^{(-1)}$ be its outward shift, and consider the face $\tau = \Delta \cap \tau'^{(-1)}$
of $\Delta$, whose integral length we denote by $r$. Then
\[ \lsize_\Sigma(\Delta^{(1)}) = s \qquad \text{and} \qquad \lsize_\Sigma(\Delta) = \max \{r, s + 3\}. \]
\end{lemma}

\noindent \emph{Remark.} The face $\tau$ is either a vertex or an edge. In the former case, its integral
length is understood to be $0$.\\

\noindent \textsc{Proof.} 
The fact that $\lsize_\Sigma(\Delta^{(1)}) = s$ follows immediately from the considerations preceding Lemma~\ref{latticesizesquare}.
As for $\lsize_\Sigma(\Delta)$, in case $r \geq s + 3$ the statement follows from Lemma~\ref{latticesize_auxiliary1}.
So assume that $r \leq s + 3$ (we reinclude the case $r = s + 3$ for the sake of the 
symmetry of the argument below). Without loss of generality we may suppose 
that $\tau' = \conv \{(0,0), (s,0) \}$. We claim that we can moreover
assume that $\tau \subset \conv \{(-1,-1), (s+2,-1) \}$, while still
keeping $\Delta^{(1)} \subset s\Sigma$. 

Assuming the claim, we can make the following reasoning.
\begin{itemize}
 \item Clearly $\Delta$ is contained in the half-plane $Y \geq -1$.
 \item Suppose that $\Delta$ contains a lattice point $(a,b)$ for which $a < - 1$.
Because $b = -1$ contradicts our claim, while $b = 0$ 
contradicts that $\Delta^{(1)} \subset s\Sigma$ (indeed, it implies that $(-1,0) \in \Delta^{(1)}$),
we must have $b \geq 1$. 
Along with the fact that $\Delta^{(1)}$ is two-dimensional (so that it must contain a lattice point on or above the
line $Y = 1$) this implies that $(0,1) \in \Delta^{(1)}$.
But then, apart from the point $(a,b)$ itself, 
all lattice points which are contained in the triangle spanned 
by $(a,b)$, $(0,0)$ and $(0,1)$ must be elements of $\Delta^{(1)}$.
The volume of this triangle being at least $1$, Pick's theorem implies that it must contain a lattice point
different from $(a,b)$, $(0,0)$ and $(0,1)$. This contradicts $\Delta^{(1)} \subset s\Sigma$.

We conclude that $\Delta$ is contained in the half-plane $X \geq -1$.

\item By applying the unimodular transformation $(i,j) \mapsto (s - i - j,j)$, one sees that the foregoing reasoning also allows to conclude that  $\Delta$ is contained in the half-plane $X + Y \leq s + 1$.
\end{itemize}
So the claim implies that $\Delta \subset (-1,-1) + (s+3)\Sigma$, and hence that $\lsize_\Sigma(\Delta) \leq s + 3$, which
together with (\ref{latticesize_trivialbound}) proves the lemma.

To prove the claim, note that because $r \leq s + 3$, again using the transformation $(i,j) \mapsto (s - i - j, j)$ if needed,
we can assume that $\tau$ is contained in the half-plane $X \geq -1$.
Let $(a,-1)$ be the right-most vertex of $\tau$. As long as $a > s+2$, we can apply a unimodular transformation
of the form $(i,j) \mapsto (i+j,j)$ to $\Delta$, while
\begin{itemize}
 \item keeping $\tau$ in the half-plane $X \geq -1$ (here we again used that $r \leq s + 3$); 
 \item keeping $\Delta^{(1)}$ inside $s\Sigma$: indeed, because $a > s + 2$ and $(s+1,0) \notin \Delta^{(1)}$, the edge of
 $\Delta$ that is right-adjacent to $\tau$ must have a slope that is smaller than $1/2$ (in absolute value), and hence the same
 must be true for the edge of $\Delta^{(1)}$ that is right-adjacent to $\tau'$.
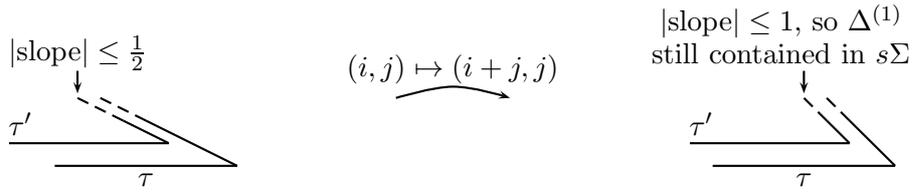
\begin{center}
\psset{unit=0.3cm}
\begin{pspicture}(-2, -2)(10,4.5)
\psline(2,-1)(10,-1)
\psline(10,-1)(6,1)
\psline[linestyle=dashed](6,1)(4,2)
\psline(0,0)(7,0)
\psline(7,0)(5,1)
\psline[linestyle=dashed](5,1)(3,2)
\rput(6,-1.6){\small $\tau$}
\rput(0.5,0.85){\small $\tau'$}
\psline{->}(3,3.2)(3,2.2)
\rput(3,4){\small $|\text{slope}| \leq \frac{1}{2}$}
\end{pspicture}
\qquad \qquad
\begin{pspicture}(-1,-2)(5,4.5)
\pscurve{->}(0,2)(2.5,2.5)(5,2)
\rput(2.5,3.3){\small $(i,j) \mapsto (i + j,j)$}
\end{pspicture}
\qquad \qquad \ \
\psset{unit=0.3cm}
\begin{pspicture}(0,0)(5,5)
\begin{pspicture}(-1, -2)(9,4.5)
\psline(1,-1)(9,-1)
\psline(9,-1)(7,1)
\psline[linestyle=dashed](7,1)(6,2)
\psline(0,0)(7,0)
\psline(7,0)(6,1)
\psline[linestyle=dashed](6,1)(5,2)
\rput(5,-1.6){\small $\tau$}
\rput(0.5,0.85){\small $\tau'$}
\psline{->}(5,3.2)(5,2.2)
\rput(4,5.4){\small $|\text{slope}| \leq 1$, so $\Delta^{(1)}$}
\rput(4,4){\small still contained in $s\Sigma$}
\end{pspicture}\end{pspicture}
\end{center}
\end{itemize}
\noindent This decreases the value of $a$ by $1$. So the claim follows by repeating this step until $a \leq s+2$.
\hfill $\blacksquare$\\

We are now ready to state and prove our recursive expression. 

\begin{theorem} \label{latticesize_classification}
 Let $\Delta$ be a two-dimensional lattice polygon. Then 
 \[ \lsize_\Sigma(\Delta) = \lsize_\Sigma(\Delta^{(1)}) + 3, \] except
 in the following situations:
 \begin{itemize}
  \item $\Delta$ is equivalent to a Lawrence prism 
  \begin{center}
  \psset{unit=0.25cm}
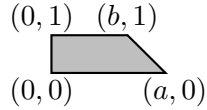
  \begin{pspicture}(-1,-1)(7,2.8)
  \pspolygon[fillstyle=solid,fillcolor=lightgray](0,0)(6,0)(4,2)(0,2)
  \rput(-0.5,-1){\small $(0,0)$}
  \rput(6.5,-1){\small $(a,0)$}  
  \rput(-0.5,2.9){\small $(0,1)$}  
  \rput(4,2.9){\small $(b,1)$}  
  \end{pspicture}
  \end{center}
  where $a=b=1$ or $2 \leq a \geq b \geq 0$, in which case $\lsize_\Sigma(\Delta^{(1)}) = -2$ and
  \[ \begin{cases} 
    \phantom{i} \lsize_\Sigma(\Delta) = a + 1 & \text{if $a = b$,} \\
    \phantom{i} \lsize_\Sigma(\Delta) = a & \text{if $a > b$;} \\
  \end{cases} \]
  \item $\Delta$ is equivalent to 
    \begin{center}
  \psset{unit=0.25cm}
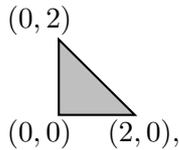
  \begin{pspicture}(-1,-1)(4,4.5)
  \pspolygon[fillstyle=solid,fillcolor=lightgray](0,0)(4,0)(0,4)
  \rput(-1,-1){\small $(0,0)$}  
  \rput(4.5,-1){\small $(2,0),$}  
  \rput(-1,5){\small $(0,2)$}  
  \end{pspicture}
  \end{center}
  in which case $\lsize_\Sigma(\Delta^{(1)}) = -2$ and $\lsize_\Sigma(\Delta) = 2$;
\item $\Delta$ is equivalent to
  \begin{center}  
 \psset{unit=0.25cm}
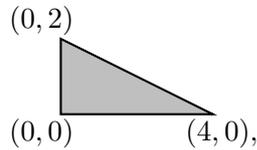
  \begin{pspicture}(-1,-1)(8,4.5)
  \pspolygon[fillstyle=solid,fillcolor=lightgray](0,0)(8,0)(0,4)
  \rput(-1,-1){\small $(0,0)$}  
  \rput(8.5,-1){\small $(4,0),$}  
  \rput(-1,5){\small $(0,2)$}  
  \end{pspicture}
  \end{center}
  in which case $\lsize_\Sigma(\Delta^{(1)}) = 0$ and $\lsize_\Sigma(\Delta) = 4$;
  \item $\Delta \simeq \square_{a,b}$ for certain $a,b \geq 2$, in which case 
  \[ \lsize_\Sigma(\Delta^{(1)}) = a + b - 4 \qquad \text{and} \qquad
  \lsize_\Sigma(\Delta) = a + b; \]
  \item there exist parallel edges
  $\tau \subset \Delta$ and $\tau' \subset \Delta^{(1)}$ whose
  supporting lines are at integral distance $1$ of each other, such that
  \[ \sharp (\tau \cap \mathbb{Z}^2) - \sharp (\tau' \cap \mathbb{Z}^2) \geq 4;\]
  in this case $\lsize_\Sigma(\Delta^{(1)}) = \sharp (\tau' \cap \mathbb{Z}^2)$ and $\lsize_\Sigma(\Delta) = \sharp (\tau \cap \mathbb{Z}^2)$.\\
  \end{itemize}
\end{theorem}

\noindent \emph{Remark.} The third case $\conv \{ (0,0), (4,0), (0,2) \}$ can 
in some sense be viewed as a special case
of the last item, with $\tau'$ having length $0$.\\

\noindent \textsc{Proof.} 
For the Lawrence prisms and the two explicit polygons 
the statement is immediate, while the polygons $\square_{a,b}$ are covered by Lemma~\ref{latticesizesquare}
and the observation that $(\square_{a,b})^{(1)} \simeq \square_{a-2,b-2}$.
The last statement follows from Lemma~\ref{latticesize_auxiliary1}.

By (\ref{latticesize_trivialbound}) it remains to show that in all other situations $\lsize_\Sigma(\Delta^{(1)}) \geq \lsize_\Sigma(\Delta) - 3$.
The cases where $\Delta^{(1)}$ is not two-dimensional can be analyzed explicitly using
Koelman's classification: see \cite[Thm.~10]{movingout} or \cite[Ch.~4]{Koelman}.
We can therefore assume
that $\Delta^{(1)}$ is two-dimensional.
Let $s = \lsize_\Sigma(\Delta^{(1)})$, so that we can suppose that
$\Delta^{(1)} \subset s \Sigma$.
If 
\begin{equation} \label{movedoutedges}
\Delta^{(1)(-1)} \subset (s\Sigma)^{(-1)}
\end{equation}
then the theorem follows because $\Delta \subset \Delta^{(1)(-1)}$ and $(s\Sigma)^{(-1)} \simeq (s + 3)\Sigma$. So let us assume that 
(\ref{movedoutedges}) is not satisfied. Without loss of generality we may then suppose that
$\Delta^{(1)(-1)}$ is not contained in the half-plane
\[ X + Y \leq s + 1.\]
This means that the edge of $s\Sigma$ connecting $(s,0)$ and $(0,s)$ cannot contain two vertices of $\Delta^{(1)}$.
But it must contain at least one vertex $v$ of $\Delta^{(1)}$: if not, $\Delta^{(1)}$ would
be contained in $(s-1)\Sigma$, contradicting $s = \lsize_\Sigma(\Delta^{(1)})$. 

 Write $v = (a, s-a)$ for some $a \in \{0, \dots, s\}$. We distinguish between two cases.
\begin{itemize}
 \item Assume that $v$ lies in the interior of the edge of $s\Sigma$ that
 connects $(s,0)$ and $(0,s)$, i.e.\ $a \notin \{0,s\}$. Let $v_1 = (a_1,b_1)$ and $v_2 = (a_2,b_2)$ be the 
 vertices of $\Delta^{(1)}$ that are adjacent to $v$, ordered counterclockwise, and for $i=1,2$ 
 let $\tau_i$ be the edge connecting $v_i$ and $v$.
 Note that $b_1 < s - a$: otherwise 
 $\Delta^{(1)}$ would be contained in $\conv \{ (0, s-a), (a, s-a), (0,s) \} \simeq a \Sigma$, which
 would contradict $s = \lsize_\Sigma(\Delta^{(1)})$. 
 This means that the outward shift $\tau_1^{(-1)}$
 must intersect the line segment spanned by $v = (a, s-a)$ and $v'= (a + 1, s-a)$. 
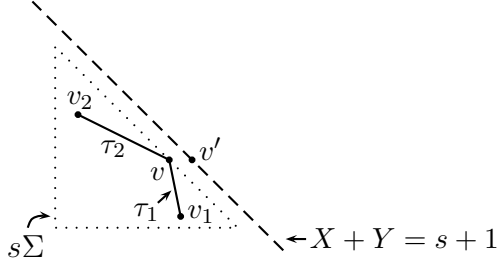
\begin{center}
\psset{unit=0.3cm}
\begin{pspicture}(-1, -1)(16,10)
\pspolygon[linestyle=dotted](0,0)(8,0)(0,8)
\psline[linestyle=dashed](10,-1)(-1,10)
\psline{->}(11,-0.5)(10.1,-0.5)
\pscircle[fillstyle=solid,fillcolor=black](5,3){0.15}
\pscircle[fillstyle=solid,fillcolor=black](6,3){0.15}
\rput(6.8,3.5){\small $v'$}
\rput(4.5,2.4){\small $v$}
\psline(5,3)(5.5,0.5)
\pscircle[fillstyle=solid,fillcolor=black](5.5,0.5){0.15}
\psline(5,3)(1,5)
\pscircle[fillstyle=solid,fillcolor=black](1,5){0.15}
\rput(6.4,0.6){\small $v_1$}
\rput(1.1,5.65){\small $v_2$}
\rput(15.3,-0.5){\small $X + Y = s + 1$}
\rput(2.55,3.5){\small $\tau_2$}
\rput(4,0.7){\small $\tau_1$}
\psline{->}(4.5,1.1)(5.1,1.65)
\rput(-1.3,-0.8){$s\Sigma$}
\pscurve{->}(-1.2,0)(-1,0.4)(-0.2,0.5)
\end{pspicture}
\end{center}
 But then 
 $b_2 \leq s - a$, otherwise $\tau_2^{(-1)}$ would
 also pass in between $v$ and $v'$, implying that $\tau_1^{(-1)}$ and $\tau_2^{(-1)}$ 
 intersect in the half-plane $X + Y \leq s + 1$: a contradiction. We conclude that $\Delta^{(1)}$ must be contained
 below the line $Y = s-a$. By symmetry of arguments, it must also lie to the left of
 $X = a$. Thus $\Delta^{(1)}$ is contained in the rectangle
 \[ \text{conv} \left\{ (0,0), (a,0), (a, s-a), (0, s-a) \right\}.\]
 Now if any of these four vertices would not appear as an actual vertex of $\Delta^{(1)}$ then
 we would again contradict $s = \lsize_\Sigma(\Delta^{(1)})$. Thus $\Delta^{(1)}$ must be exactly this rectangle,
 and $\Delta^{(1)(-1)} \simeq \square_{a+2,s-a+2}$. 
 The case
 $\Delta = \Delta^{(1)(-1)}$ being among our exceptions,
we can assume that at least
 one of the four vertices of $\Delta^{(1)(-1)}$ does not appear as an actual vertex of $\Delta$. But then 
 $\lsize_\Sigma(\Delta) \leq s + 3$, as desired. 
 
 \item Assume that $v$ is an endpoint of the edge of $s\Sigma$ connecting $(s,0)$ and $(0,s)$, i.e.\ $a \in \{0,s\}$.
 Without loss of generality we may assume that $a = s$. Again let $v_1 = (a_1,b_1)$ and $v_2 = (a_2,b_2)$ be the 
 vertices of $\Delta^{(1)}$ that are adjacent to $v$, ordered counterclockwise, and for $i=1,2$ 
 let $\tau_i$ be the edge connecting $v_i$ and $v$. 
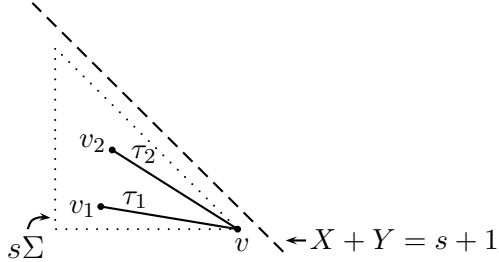
\begin{center}
\psset{unit=0.3cm}
\begin{pspicture}(-1, -1)(16,10)
\pspolygon[linestyle=dotted](0,0)(8,0)(0,8)
\psline[linestyle=dashed](10,-1)(-1,10)
\psline{->}(11,-0.5)(10.1,-0.5)
\pscircle[fillstyle=solid,fillcolor=black](8,0){0.15}
\rput(8.2,-0.7){\small $v$}
\rput(15.3,-0.5){\small $X + Y = s + 1$}
\psline(8,0)(2.5,3.5)
\psline(8,0)(2,1)
\pscircle[fillstyle=solid,fillcolor=black](2.5,3.5){0.15}
\pscircle[fillstyle=solid,fillcolor=black](2,1){0.15}
\rput(1.65,3.8){\small $v_2$}
\rput(1.2,0.9){\small $v_1$}
\rput(3.5,1.3){\small $\tau_1$}
\rput(3.9,3.32){\small $\tau_2$}
\rput(-1.3,-0.8){$s\Sigma$}
\pscurve{->}(-1.2,0)(-1,0.4)(-0.2,0.5)
\end{pspicture}
\end{center}
 We claim that $v_1 = (0,0)$, i.e.\ $a_1 = b_1 = 0$. Indeed:
  \begin{itemize}
   \item Assume that $b_1 = 0$. Then $\tau_1^{(-1)}$ is the line $Y = -1$. Since $\tau_2^{(-1)}$ must intersect
   this line in a lattice point outside the half-plane $X + Y \leq s + 1$ we find (as in the proof of Lemma~\ref{latticesize_auxiliary2}) that
   $\tau_2$ has slope at most $1/2$ (in absolute value), i.e.\ $a_2 \leq s - 2b_2$.
   From this it follows that $a_1 = 0$: if not, the unimodular transformation
   $(i,j) \mapsto (i + j-1, j)$ maps $\Delta^{(1)}$ inside $(s-1)\Sigma$, contradicting $s = \lsize_\Sigma(\Delta^{(1)})$.
   \item Assume that $b_1 \neq 0$. If $a_2 \leq s - 2b_2$ then we would
   again find a contradiction with $s = \lsize_\Sigma(\Delta^{(1)})$. 
   Therefore $a_2 > s - 2b_2$, and by symmetry of arguments also $a_1 < s - 2b_1$. But then
   $\tau_1^{(-1)}$ passes through or above the point $(s + 2, -1)$, while
   $\tau_2^{(-1)}$ passes through or to the left of $(s + 2, -1)$.  
   Taking into account their respective slopes, one sees that these lines 
   must intersect in the half-plane
   $X + Y \leq s + 1$: a contradiction. 
   So this case cannot occur.
 \end{itemize}
 Thus $\tau_1 = \conv \{ (0,0), (s,0) \}$. Now consider the face $\tau = \tau_1^{(-1)} \cap \Delta$ of $\Delta$. The
 case $\sharp (\tau \cap \mathbb{Z}^2) \geq s + 4$ being among our exceptions, we can assume
 that $\sharp (\tau \cap \mathbb{Z}^2) \leq s + 3$. The theorem then follows from 
 Lemma~\ref{latticesize_auxiliary2}. \hfill $\blacksquare$\\
\end{itemize}

Theorem~\ref{latticesize_classification} 
gives a recursive method for computing the lattice size with respect to $\Sigma$ in practice.
For example, let $\Delta$ be the lattice polygon below.
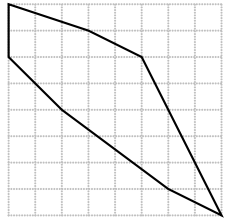
\begin{center}
\psset{unit=0.35cm}
\begin{pspicture}(0,-6)(8,2)
\psgrid[subgriddiv=1,griddots=10,gridcolor=lightgray,gridlabels=0pt](0,-6)(8,2)
\pspolygon(0,0)(2,-2)(6,-5)(8,-6)(5,0)(3,1)(0,2)
\end{pspicture}
\end{center}
By taking consecutive interiors, we find the following `onion skins'.
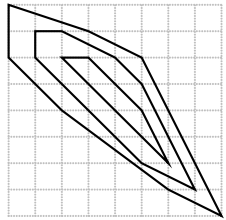
\begin{center}
\psset{unit=0.35cm}
\begin{pspicture}(0,-6)(8,2)
\psgrid[subgriddiv=1,griddots=10,gridcolor=lightgray,gridlabels=0pt](0,-6)(8,2)
\pspolygon(0,0)(2,-2)(6,-5)(8,-6)(5,0)(3,1)(0,2)
\pspolygon(1,0)(5,-4)(7,-5)(5,-1)(4,0)(2,1)(1,1)
\pspolygon(2,0)(6,-4)(5,-2)(3,0)
\end{pspicture}
\end{center}
The inner polygon is (equivalent to) a Lawrence prism with $a = 4$ and $b = 2$, while
the subsequent steps are not exceptional. We find $\lsize_\Sigma(\Delta) = \lsize_\Sigma(\emptyset) + 6 + 3 + 3 = 10$.
We remark that this is in fact a rephrasing of Schicho's algorithm for simplifying rational parametrizations
of toric surfaces \cite[\S4]{schicho}. Whereas Schicho proved that the output of the algorithm is at worst twice the lattice size \cite[Thm.\,10]{schicho}, our result
shows that the result is actually optimal.\\

\noindent A Magma implementation of this method can be found in the file \verb"basic_commands.m"
that accompanies \cite{linearpencils}. For instance, the above example can be treated as follows:

\begin{verbatim}
> load "basic_commands.m";
Loading "basic_commands.m"
> P := LatticePolytope([<8,0>,<6,1>,<2,4>,<0,6>,<0,8>,<3,7>,<5,6>]);
> LatticeSizeRecursiveSigma(P);
10
\end{verbatim}

\section{A recursive formula for $\lsize_\square(\Delta)$}  \label{section_recursivesquare}

Some basic properties of the lattice size with respect to $\square$ are that
\[ \lsize_\square(\conv\{ (0,0), (d,0) \}) = d \]
for any $d \in \mathbb{Z}_{\geq 0}$ (in particular every lattice polygon that contains a line segment of integral length $d$ must have
lattice size at least $d$ with respect to $\square$), and that for each
non-empty lattice polygon $\Delta$ we have
\begin{equation} \label{latticesizechain}
 \lwidth(\Delta) \leq \lsize_\square(\Delta) \leq \lsize_\Sigma(\Delta) \leq 2 \lsize_\square (\Delta).
\end{equation}
By Lemma~\ref{spectrum} the first two inequalities become equalities for (and only for) $\Delta \simeq d\Sigma$ with $d \in \mathbb{Z}_{\geq 0}$.

The aim is again to relate $\lsize_\square(\Delta)$ to $\lsize_\square(\Delta^{(1)})$. Our treatment is very similar to that of the previous section. Because 
$(d \square)^{(1)} \simeq (d-2)\square$ for $d \geq 2$, we have that
\begin{equation} \label{latticesizesquare_trivialbound}
 \lsize_\square(\Delta^{(1)}) \leq \lsize_\square(\Delta) - 2
\end{equation}
as soon as $\Delta$ is two-dimensional (this includes the case where $\Delta^{(1)} = \emptyset$, which
can be verified explicitly). Typically one expects equality to hold, so our task amounts to classifying the exceptions. We again rely on two
auxiliary lemmas. The first is a literal rephrasing of Lemma~\ref{latticesize_auxiliary1}:
\begin{lemma} \label{latticesizesquare_auxiliary1}
Assume that there exist parallel edges $\tau \subset \Delta$ and $\tau' \subset \Delta^{(1)}$ whose supporting
lines are at integral distance $1$ of each other, of respective lengths $r$ and $s$. If $r \geq s + 3$
then $\lsize_\square(\Delta^{(1)}) = s$ and $\lsize_\square(\Delta) = r$.
\end{lemma}
\noindent \textsc{Proof.} By Lemma~\ref{latticesize_auxiliary1} we know that $\lsize_\Sigma(\Delta^{(1)}) = s$ and $\lsize_\Sigma(\Delta) = r$, so by
(\ref{latticesizechain}) we find $\lsize_\square(\Delta^{(1)}) \leq s$ and $\lsize_\square(\Delta) \leq r$. Equality follows from the considerations at the beginning of this section.
\hfill $\blacksquare$\\

\noindent (Instead of invoking Lemma~\ref{latticesize_auxiliary1} one can also just copy its proof, basically.) Our second lemma is analogous to Lemma~\ref{latticesize_auxiliary2}, but the statement is slightly more subtle:

\begin{lemma} \label{latticesizesquare_auxiliary2}
Assume that $\Delta^{(1)}$ is two-dimensional. Let $s \geq 1$ be an integer such that $\Delta^{(1)} \subset s \square$, and assume that
$\Delta^{(1)}$ has at least one edge in common with $s\square$. Choose such an edge $\tau'$ for which the integral length $r$ of the face $\tau = \Delta \cap \tau'^{(-1)}$
of $\Delta$ is maximal. Then
\[ \lsize_\square(\Delta^{(1)}) = s \qquad \text{and} \qquad \lsize_\square(\Delta) = \max \{r, s + 2\}. \] 
\end{lemma}

\noindent \textsc{Proof.} The fact that $\lsize_\square(\Delta^{(1)}) = s$ follows immediately from the considerations
at the beginning of this section. As for $\lsize_\square(\Delta)$, in case $r \geq s + 3$ the statement follows from
Lemma~\ref{latticesizesquare_auxiliary1}. So assume that $r \leq s + 2$. Without loss of generality
we may suppose that $\tau' = \conv \{ (0,0), (s,0) \}$. In complete analogy with the proof of Lemma~\ref{latticesize_auxiliary2} we can moreover
assume that $\tau \subset \conv \{ (-1,-1), (s+1,-1) \}$, while still keeping $\Delta^{(1)} \subset s \square$. 
Still copying the reasoning from that proof, we conclude that $\Delta$
must be in the half-planes $Y \geq -1$, $X \geq -1$ and $X \leq s + 1$. 
\begin{center}
\psset{unit=0.35cm}
\begin{pspicture}(-2, -2)(8,12.5)
\pspolygon[linestyle=dotted](0,0)(6,0)(6,6)(0,6)
\psline[linestyle=dotted](-5,-1)(8,-1)
\psline[linestyle=dotted](-1,-2)(-1,11)
\psline[linestyle=dotted](7,-2)(7,11)
\psline(0,0)(6,0)
\psline(1,-1)(6,-1)
\rput(-3.7,-0.3){\small $Y \geq -1$}
\rput{90}(-0.2,9){\small $X \geq -1$}
\rput{90}(6.3,9.5){\small $X \leq s+1$}
\rput(3.4,0.8){\small $\tau'$}
\rput(3.6,-1.6){\small $\tau$}
\rput(3.5,5.1){\small $s\square$}
\end{pspicture}
\end{center}
Now suppose that $\Delta$ contains a lattice point
$(a,b)$ for which $b > s + 1$. If $0 \leq a \leq s$ then the point $(a,s+1)$ is contained in the triangle spanned by $(a,b)$, $(0,0)$ and $(s,0)$, therefore it must be
contained in $\Delta^{(1)}$, contradicting that $\Delta^{(1)} \subset s\square$. We can therefore make the following case distinction:
\begin{itemize}
  \item $\Delta$ is contained in the half-plane $Y \leq s + 1$. But 
  this means that $\Delta \subset (-1,-1) + (s+2)\square$ and hence that $\lsize_\square(\Delta) \leq s + 2$, which together
with (\ref{latticesizesquare_trivialbound}) allows us to conclude.
  \item $\Delta$ contains a point $(-1,b)$ with $b > s+1$. By considering the convex hull of $(-1,b)$, $(0,0)$ and $(s,0)$ this implies
  that $s\Sigma \subset \Delta^{(1)}$. Now if 
  \begin{itemize}
    \item the latter inclusion would be strict, or 
    \item if $b > s + 2$, 
   \end{itemize}
   then one would obtain that $(0,s+1) \in \Delta^{(1)} \subset s\square$, a contradiction.
  Therefore $b = s + 2$ and $\Delta^{(1)} = s\Sigma$. In particular $\Delta^{(1)}$ also
  has the vertical edge $\conv \{ (0,0), (0,s) \}$ in common with $s\square$. This means that $(-1,-1) \notin \Delta$, for otherwise
  the corresponding face of $\Delta$ would contain $\conv \{ (-1,-1), (-1,s+2) \}$ which has integral length $s + 3$, contradicting the maximality of $r$.
  But then the unimodular transformation $(i,j) \mapsto (i,i+j)$ maps $\Delta$ inside $(-1,-1) + (s+2)\square$. Hence $\lsize_\square(\Delta) \leq s + 2$, which together
with (\ref{latticesizesquare_trivialbound}) allows us to conclude.
  \item $\Delta$ contains a point $(s+1,b)$ with $b > s+1$. This case follows from the previous one, by symmetry.
\end{itemize}
This proves the lemma.
\hfill $\blacksquare$\\

\noindent In the statement of Lemma~\ref{latticesizesquare_auxiliary2}, the condition of maximality is necessary. For instance let $\Delta$
be the polygon
  \begin{center}
  \psset{unit=0.25cm}
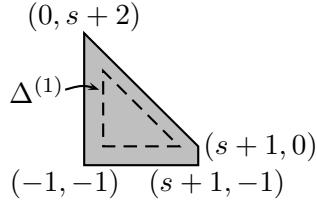
  \begin{pspicture}(-2,-2.2)(6,7.2)
  \pspolygon[fillstyle=solid,fillcolor=lightgray](-1,-1)(5,-1)(5,0)(-1,6)
  \pspolygon[linestyle=dashed](0,0)(4,0)(0,4)
  \rput(-2,-2){\small $(-1,-1)$}
  \rput(6,-2){\small $(s+1,-1)$}
  \rput(8.3,0){\small $(s+1,0)$}
  \rput(-1,6.9){\small $(0,s + 2)$}
  \pscurve{->}(-2,3)(-1.1,3.2)(-0.2,3)
  \rput(-3.5,3){\small $\Delta^{(1)}$}
  \end{pspicture}
  \end{center}
  so that $\Delta^{(1)} = s\Sigma \subset s\square$. Both 
  $\conv \{(0,0), (0,s) \}$ and $\conv \{ (0,0), (s,0) \}$ are common edges, 
  but the corresponding faces $\tau$ of $\Delta$ have different integral lengths, namely $s + 3$ resp.\ $s+2$. 
  So in this case the lattice size of $\Delta$ with respect to $\square$ is $s + 3$.
  
Let us include the following corollary to (the proof of) Lemma~\ref{latticesizesquare_auxiliary2}, for use in Section~\ref{section_minproductorder}. 
Define a \emph{horizontal} resp.\ \emph{vertical skewing} as a unimodular transformation of the form
\[ \begin{pmatrix} i \\ j \\ \end{pmatrix} \mapsto \begin{pmatrix} \pm 1 & a \\ 0 & 1 \\ \end{pmatrix} \begin{pmatrix} i \\ j \\ \end{pmatrix} + \begin{pmatrix} b \\ 0 \\ \end{pmatrix} \qquad \text{resp.} 
 \qquad  \begin{pmatrix} i \\ j \\ \end{pmatrix} \mapsto \begin{pmatrix} 1 & 0 \\ a & \pm 1 \\ \end{pmatrix} \begin{pmatrix} i \\ j \\ \end{pmatrix} + \begin{pmatrix} 0 \\ b \\ \end{pmatrix}    \]
 for some $a,b \in \mathbb{Z}$ (i.e.\ leaving the second resp.\ first coordinate invariant).

\begin{corollary} \label{skewingcorollary}
 Assume that $\Delta^{(1)}$ is two-dimensional and
 contained in $\lsize_\square(\Delta^{(1)}) \cdot \square$. Suppose that these polygons have a unique edge in common.  
 If this edge is horizontal (resp.\ vertical) then there exists a horizontal (resp.\ vertical) skewing
 $\varphi$ for which $\varphi(\Delta) \subset (-1,-1) + \lsize_\square(\Delta) \cdot \square$.
\end{corollary}

\noindent \textsc{Proof.} Let $\tau'$ be the common edge with $s \square$ where $s = \lsize_\square(\Delta^{(1)})$, and let $\tau = \tau'^{(-1)} \cap \Delta$.
Denote the integral length of the latter by $r$. Without loss of generality we can assume that $\tau'$ is a horizontal edge of $s\square$. 

We actually claim the stronger statement that there exists a horizontal skewing $\varphi$ such that
\[ \varphi(\Delta) \subset [-1,\lsize_\square(\Delta) - 1] \times [-1, s + 1] .\]
To prove this it suffices to assume that $\tau'$ is the bottom edge of $s\square$, so that we are in the set-up from the proof of Lemma~\ref{latticesizesquare_auxiliary2}. We make a
case distinction:
\begin{itemize}
\item either $r \leq s + 2$, in which case the proof of Lemma~\ref{latticesizesquare_auxiliary2} yields that
$\lsize_\square(\Delta) = s + 2$ and, through the proof of Lemma~\ref{latticesize_auxiliary2}, that 
there exists a horizontal skewing $\varphi$ such that 
\[ \varphi(\Delta^{(1)}) \subset s\square \qquad \text{and} \qquad \varphi(\Delta) \subset (-1,-1) + (s+2)\square: \] 
indeed, the unicity of $\tau'$ excludes that $s\Sigma \subset \varphi(\Delta^{(1)})$, so
the last part of the proof of Lemma~\ref{latticesizesquare_auxiliary2} can be omitted; 
\item or $r \geq s + 3$, in which case Lemma~\ref{latticesizesquare_auxiliary1} yields that $\lsize_\square(\Delta) = r$ and, through
the the proof of Lemma~\ref{latticesize_auxiliary1}, that there exists a horizontal skewing $\varphi$ such that
$\varphi(\Delta) \subset (-1,-1) + r\square$.
\end{itemize}
In both cases the claim follows. \hfill $\blacksquare$\\

We now state and prove our recursive expression.

\begin{theorem} \label{classificationsquare}
Let $\Delta$ be a two-dimensional lattice polygon. Then 
\[ \lsize_\square(\Delta) = \lsize_\square(\Delta^{(1)}) + 2,\] except in the following
situations:
\begin{itemize}
  \item $\Delta$ is equivalent to a Lawrence prism 
  \begin{center}
  \psset{unit=0.25cm}
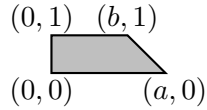
  \begin{pspicture}(-1,-1)(7,2.8)
  \pspolygon[fillstyle=solid,fillcolor=lightgray](0,0)(6,0)(4,2)(0,2)
  \rput(-0.5,-1){\small $(0,0)$}
  \rput(6.5,-1){\small $(a,0)$}  
  \rput(-0.5,2.9){\small $(0,1)$}  
  \rput(4,2.9){\small $(b,1)$}  
  \end{pspicture}
  \end{center}
  where $2 \leq a \geq b \geq 0$, in which case $\lsize_\square(\Delta^{(1)}) = -1$ and $\lsize_\square(\Delta) = a$;

  \item $\Delta$ is equivalent to
      \begin{center}
  \psset{unit=0.25cm}
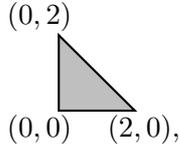
  \begin{pspicture}(-1,-1.5)(4,4.5)
  \pspolygon[fillstyle=solid,fillcolor=lightgray](0,0)(4,0)(0,4)
  \rput(-1,-1){\small $(0,0)$}  
  \rput(4.5,-1){\small $(2,0),$}  
  \rput(-1,5){\small $(0,2)$}  
  \end{pspicture}
   \end{center}
  in which case $\lsize_\square(\Delta^{(1)}) = -1$ and $\lsize_\square(\Delta) = 2$;
  \item 
   $\Delta$ is equivalent to one of
    \begin{center}
  \psset{unit=0.25cm}
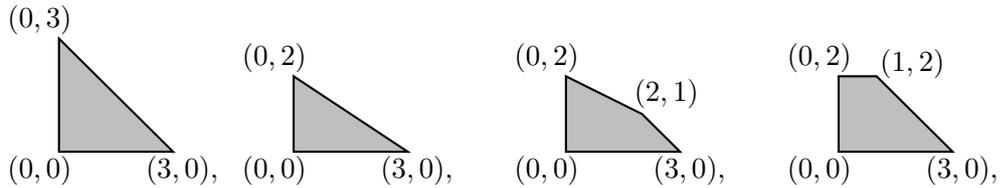
  \begin{pspicture}(-1,-1.5)(4,6.7)
  \pspolygon[fillstyle=solid,fillcolor=lightgray](0,0)(6,0)(0,6)
  \rput(-1,-1){\small $(0,0)$}  
  \rput(6.5,-1){\small $(3,0),$}  
  \rput(-1,7){\small $(0,3)$}  
  \end{pspicture}
   \qquad \qquad
  \begin{pspicture}(-1,-1.5)(6,4.5)
  \pspolygon[fillstyle=solid,fillcolor=lightgray](0,0)(6,0)(0,4)
  \rput(-1,-1){\small $(0,0)$}  
  \rput(6.5,-1){\small $(3,0),$}  
  \rput(-1,5){\small $(0,2)$}  
  \end{pspicture}
   \qquad \qquad
    \begin{pspicture}(-1,-1.5)(6,4.5)
  \pspolygon[fillstyle=solid,fillcolor=lightgray](0,0)(6,0)(4,2)(0,4)
  \rput(-1,-1){\small $(0,0)$}  
  \rput(6.5,-1){\small $(3,0),$}  
  \rput(-1,5){\small $(0,2)$}  
    \rput(5.3,3.1){\small $(2,1)$}
  \end{pspicture}
   \qquad \qquad
    \begin{pspicture}(-1,-1.5)(6,4.5)
  \pspolygon[fillstyle=solid,fillcolor=lightgray](0,0)(6,0)(2,4)(0,4)
  \rput(-1,-1){\small $(0,0)$}  
  \rput(6.5,-1){\small $(3,0),$}  
  \rput(-1,5){\small $(0,2)$}  
  \rput(3.9,4.7){\small $(1,2)$}
  \end{pspicture}
  \end{center}
   in which case $\lsize_\square(\Delta^{(1)}) = 0$ and $\lsize_\square(\Delta) = 3$;
  \item 
   $\Delta$ is equivalent to 
    \begin{center}
  \psset{unit=0.25cm}
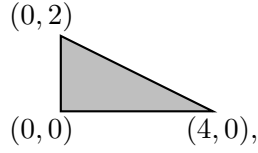
    \begin{pspicture}(-1,-1)(8,4.5)
  \pspolygon[fillstyle=solid,fillcolor=lightgray](0,0)(8,0)(0,4)
  \rput(-1,-1){\small $(0,0)$}  
  \rput(8.5,-1){\small $(4,0),$}  
  \rput(-1,5){\small $(0,2)$}  
  \end{pspicture}
  \end{center}
  in which case $\lsize_\square(\Delta^{(1)}) = 0$ and $\lsize_\square(\Delta) = 4$;
  \item there exist parallel edges
  $\tau \subset \Delta$ and $\tau' \subset \Delta^{(1)}$ whose
  supporting lines are at integral distance $1$ of each other, such that
  \[ \sharp (\tau \cap \mathbb{Z}^2) - \sharp (\tau' \cap \mathbb{Z}^2) \geq 3;\]
  in this case $\lsize_\square(\Delta^{(1)}) = \sharp (\tau' \cap \mathbb{Z}^2)$ and $\lsize_\square(\Delta) = \sharp (\tau \cap \mathbb{Z}^2)$.\\
  \end{itemize}
\end{theorem}

\noindent \emph{Remark.} Except for $2\Sigma$, the explicitly given polygons can in some sense be viewed as special cases
of the last item, with $\tau'$ having length $0$.\\

\noindent \textsc{Proof.} 
For the Lawrence prisms and the six explicitly given polygons, the theorem is immediate, while
the last statement follows from Lemma~\ref{latticesizesquare_auxiliary1}.

By (\ref{latticesizesquare_trivialbound}) it remains to show that in all other situations $\lsize_\square(\Delta^{(1)}) \geq \lsize_\square(\Delta) - 2$.
The cases where $\Delta^{(1)}$ is not two-dimensional can again be analyzed explicitly using Koelman's classification: see \cite[Thm.~10]{movingout} or \cite[Ch.~4]{Koelman}. We can
therefore assume that $\Delta^{(1)}$ is two-dimensional. Let $s = \lsize_\square(\Delta^{(1)})$ and suppose that $\Delta^{(1)} \subset s \square$.
If 
\begin{equation} \label{movedoutedgessquare}
\Delta^{(1)(-1)} \subset (s\square)^{(-1)}
\end{equation}
then the theorem follows because $\Delta \subset \Delta^{(1)(-1)}$ and $(s\square)^{(-1)} \simeq (s + 2)\square$. So let us assume that 
(\ref{movedoutedgessquare}) is not satisfied. Without loss of generality we may then suppose that
$\Delta^{(1)(-1)}$ is not contained in the half-plane $X \leq s + 1$.
By using a translation if needed, we can assume that the both the lower edge and the right edge of $s\square$ contain
at least one vertex of $\Delta^{(1)}$. 
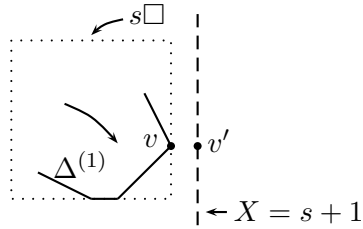
\begin{center}
\psset{unit=0.35cm}
\begin{pspicture}(-1, -1)(7,7)
\pspolygon[linestyle=dotted](0,0)(6,0)(6,6)(0,6)
\psline[linestyle=dashed](7,-1)(7,7)
\psline(1,1)(3,0)
\psline(3,0)(4,0)
\psline(4,0)(6,2)
\psline(6,2)(5,4)
\rput(10.8,-0.52){\small $X = s + 1$}
\pscircle[fillstyle=solid,fillcolor=black](6,2){0.15}
\pscircle[fillstyle=solid,fillcolor=black](7,2){0.15}
\psline{->}(8.1,-0.5)(7.3,-0.5)
\pscurve{->}(2,3.6)(3,3.1)(4,2.1)
\rput(5.25,2.2){\small $v$}
\rput(7.8,2.2){\small $v'$}
\rput(5.1,7){\small $s\square$}
\rput(2.6,1.3){\small $\Delta^{(1)}$}
\pscurve{->}(4.2,6.8)(3.5,6.6)(3,6.2)
\end{pspicture}
\end{center}
By our assumption the right edge  
then contains \emph{exactly} one such vertex, that we denote by $v = (s,a)$,
for some $a \in \{ 0, \dots, s \}$.

We first reduce to the case where $a = 0$. Suppose that $a > 0$ and let $v_1 = (a_1,b_1)$ and $v_2 = (a_2,b_2)$ be the vertices
of $\Delta^{(1)}$ that are adjacent to $v$, ordered counterclockwise. For $i = 1,2$ let $\tau_i$ be the edge connecting $v_i$ and $v$. By
our assumption that the lower edge of $s\square$ contains at least one vertex of $\Delta^{(1)}$ we have that $b_1 < a$. This means
that the outward shift $\tau_1^{(-1)}$ must intersect the line segment spanned by $v = (s,a)$ and $v' = (s+1,a)$. But then $b_2 \leq a$, otherwise
$\tau_2^{(-1)}$ would also pass in between $v$ and $v'$, implying that $\tau_1^{(-1)}$ and $\tau_2^{(-1)}$ intersect in the half-plane $X \leq s + 1$: a contradiction.
We conclude that $\Delta^{(1)}$ lies in the half-plane $Y \leq a$. But then a vertical flip followed by a vertical translation positions $v$ at $(s,0)$, while
leaving our other assumptions unaffected.

So we can assume that $v = (s,0)$. We claim that this implies that at least one of $\conv \{ (0,0), (s,0) \}$ or 
$\conv \{ (0,s), (s,0) \}$ appears as an edge of $\Delta^{(1)}$. Assuming the claim we can conclude quickly. Indeed, in the former case we see that $\Delta^{(1)}$ has an
edge in common with $s \square$, so that the theorem follows from Lemma~\ref{latticesizesquare_auxiliary2} (using that
we excluded the cases where $r \geq s + 3$). In the latter case either $(i,j) \mapsto (i+j,j)$ or $(i,j) \mapsto (i,i+j-s)$ positions 
$\Delta^{(1)}$ inside $s\square$ in such a way that there is an edge in common: 
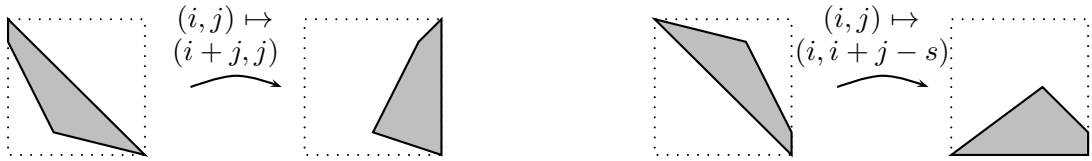
\begin{center}
\psset{unit=0.3cm}
\begin{pspicture}(-1, -1)(20,7)
\pspolygon[linestyle=dotted](0,0)(6,0)(6,6)(0,6)
\pspolygon[fillstyle=solid,fillcolor=lightgray](0,6)(6,0)(2,1)(1,3)(0,5)
\pscurve{->}(8,3)(10,3.5)(12,3)
\rput(9.5,6){\small $(i,j) \mapsto$}
\rput(9.55,4.5){\small $(i+j,j)$}
\pspolygon[linestyle=dotted](13,0)(19,0)(19,6)(13,6)
\pspolygon[fillstyle=solid,fillcolor=lightgray](19,6)(19,0)(16,1)(17,3)(18,5)
\end{pspicture}
\qquad \qquad \quad
\begin{pspicture}(-1, -1)(20,7)
\pspolygon[linestyle=dotted](0,0)(6,0)(6,6)(0,6)
\pspolygon[fillstyle=solid,fillcolor=lightgray](6,0)(0,6)(4,5)(5,3)(6,1)
\pscurve{->}(8,3)(10,3.5)(12,3)
\rput(9.5,6){\small $(i,j) \mapsto$}
\rput(9.55,4.5){\small $(i,i+j-s)$}
\pspolygon[linestyle=dotted](13,0)(19,0)(19,6)(13,6)
\pspolygon[fillstyle=solid,fillcolor=lightgray](19,0)(13,0)(17,3)(18,2)(19,1)
\end{pspicture}
\end{center}
So the theorem again follows from Lemma~\ref{latticesizesquare_auxiliary2}.

To prove the claim, as before let $v_1 = (a_1,b_1)$ and $v_2 = (a_2,b_2)$ be the vertices of $\Delta^{(1)}$ that
are adjacent to $v$, ordered counterclockwise, and denote by $\tau_1, \tau_2$ the corresponding edges. We make a case distinction.
\begin{itemize}
  \item Assume that $b_1 = 0$. Then $\tau_1^{(-1)}$ is the line $Y = -1$. Since $\tau_2^{(-1)}$ must
  intersect this line in a lattice point outside the half-plane $X \leq s + 1$ we find that $\tau_2$ has slope at most $1$ (in absolute value), i.e.\ $a_2 \leq s - b_2$.
  It follows that $\Delta^{(1)} \subset s\Sigma$. Now:
  \begin{itemize}
    \item if $(0,0) \in \Delta^{(1)}$ or $(0,s) \in \Delta^{(1)}$ then the claim follows;
    \item if not then the transformation $(i,j) \mapsto (i + j - 1,j)$ maps $\Delta^{(1)}$ inside $(s-1) \square$, contradicting that
    $s = \lsize_\square(\Delta^{(1)})$.
  \end{itemize}
  \item Assume that $b_1 \neq 0$. Then $\tau_1$ and $\tau_2$ cannot lie at opposite sides of the line connecting $(s,0)$ and $(0,s)$, i.e.\
  one cannot simultaneously have $a_1 < s - b_1$ and $a_2 > s - b_2$, because otherwise $\tau_1^{(-1)}$ and $\tau_2^{(-1)}$ would intersect 
  in the half-plane $X \leq s + 1$. But then either $a_2 \leq s - b_2$, in which case $\Delta^{(1)} \subset s\Sigma$ and we can proceed as before,
  or $a_1 \geq s - b_1$, in which case the situation is entirely analogous.
\end{itemize}
This completes the proof.
\hfill $\blacksquare$\\

Theorem~\ref{classificationsquare} gives a recursive method for computing the lattice size with respect to $\square$ in practice.
Using the example $\Delta$ from the end of the previous section, we see that $\lsize_\square(\Delta) = \lsize_\square(\emptyset) + 5 + 2 + 2 = 8$.
A Magma implementation can be found in the file \verb"basic_commands.m" that accompanies~\cite{linearpencils}. For instance, the foregoing
example can be treated as follows:

\begin{verbatim}
> load "basic_commands.m";
Loading "basic_commands.m"
> P := LatticePolytope([<8,0>,<6,1>,<2,4>,<0,6>,<0,8>,<3,7>,<5,6>]);
> LatticeSizeRecursiveSquare(P);
8
\end{verbatim}

We include an immediate corollary to the above proof, for use in the next section. Let $s \in \mathbb{Z}_{\geq 1}$. Then by a \emph{slice} of $s \square$ we mean
a line segment of the form $\conv \{ (a,0), (a,s) \}$ or $\conv \{ (0,a), (s,a) \}$ for some $a \in \{0, \dots, s\}$. By a \emph{diagonal} we
mean $\conv \{ (0,s), (s,0) \}$ or $\conv \{(0,0), (s,s) \}$. Then we have:

\begin{corollary} \label{slicingcorollary}
  Suppose that $\Delta^{(1)}$ is two-dimensional and contained in $\lsize_\square(\Delta^{(1)}) \cdot \square$.
  Assume that there is no edge of $\Delta^{(1)}$ that is a slice or a diagonal of the latter. Then
  $\lsize_\square(\Delta) = \lsize_\square(\Delta^{(1)}) + 2$ and $\Delta \subset (-1,-1) + \lsize_\square(\Delta) \cdot \Delta$.
\end{corollary}

\section{A minimum with respect to the product order} \label{section_minproductorder}

This section is devoted to our combinatorial version of Conjecture~\ref{conjecturebidegree}, namely
that for each non-empty lattice polygon $\Delta$ the set $S_{1,1}(\Delta)$ admits a minimum with respect to the product order $\leq \times \leq$ on $\mathbb{N}^2$.
It suffices to show that
$\Delta$ admits a unimodular copy inside the rectangle
\begin{equation} \label{importantrectangle}
 \square_{\lwidth(\Delta), \lsize_\square(\Delta)  } = [0, \lwidth(\Delta) ] \ \times \ [0, \lsize_\square(\Delta) ].
\end{equation}
Indeed, then $(\lwidth(\Delta), \lsize_\square(\Delta)) \in S_{1,1}(\Delta)$, and from the respective definitions of $\lwidth(\Delta)$ and $\lsize_\square(\Delta)$ it is clear that
this concerns a minimum with respect to the product order.

We need the following properties of the lattice width:
\begin{lemma} \label{lubbesschicholemma}
    If $\Delta^{(1)} \neq \emptyset$ then $\lwidth(\Delta) = \lwidth(\Delta^{(1)}) + 2$, except if
    $\Delta \simeq d\Sigma$ for some $d \geq 3$, in which case $d = \lwidth(\Delta) = \lwidth(\Delta^{(1)}) + 3$.
    If moreover $\Delta^{(1)} \not \simeq d\Sigma$ for any $d \geq 0$ and
    \[ \Delta^{(1)} \ \subset \ [0, \lwidth(\Delta^{(1)}) ] \ \times \ \mathbb{R} \]
    then
    \[ \Delta \ \subset \ [-1, \lwidth(\Delta^{(1)}) + 1 ] \ \times \ \mathbb{R} \ = \ [-1, \lwidth(\Delta) - 1 ] \ \times \ \mathbb{R}. \]
 \end{lemma}
\noindent \textsc{Proof.} See~\cite[Thm.\,13]{lubbesschicho}, where the second statement is phrased as follows: an optimal viewangle for $\Delta^{(1)}$ is
also an optimal viewangle for $\Delta$. \hfill $\blacksquare$\\

Due to the special role of standard triangles, we treat the following case separately:

\begin{lemma}
  Let $\Delta$ be a two-dimensional lattice polygon such that $\Delta^{(1)} \simeq d \Sigma$ for some $d \geq 1$. Then there exists a unimodular transformation mapping
  $\Delta$ inside (\ref{importantrectangle}).
\end{lemma}

\noindent \textsc{Proof.} If $\Delta^{(1)} = d \Sigma$ then $\Delta \subset (d\Sigma)^{(-1)} \simeq (d+3) \Sigma$.
If equality holds then $\lwidth(\Delta) = \lsize_\square(\Delta) = d + 3$ and $\Delta$ is indeed contained
    in a box of size $(d+3) \times (d+3)$. If not then at least one of the vertices of $(d+3)\Sigma$ is not contained in $\Delta$. By applying a unimodular transformation
    if needed we can assume that it concerns the right-most vertex. We now make a case distinction:
    \begin{itemize}
    \item If the left-most edge of $(d+3)\Sigma$ is contained in $\Delta$, then $\lsize_\square(\Delta) = d + 3$ by Lemma~\ref{latticesizesquare_auxiliary1}. On the other
    hand $\lwidth(\Delta) = d + 2$ by Lemma~\ref{lubbesschicholemma}. We see that $\Delta$ is contained in a box of size $(d+2) \times (d+3)$, as wanted.
    
    (Remark: the example following the proof of Lemma~\ref{latticesizesquare_auxiliary2} is of this kind.)
    \item If the left-most edge does not appear, then without loss of generality we can assume that the top vertex is missing. Then 
    $\lsize_\square(\Delta) = d + 2$ by Lemma~\ref{latticesizesquare_auxiliary2}, while still $\lwidth(\Delta) = d + 2$. We see that $\Delta$ is contained in a box of size $(d+2) \times (d+2)$, as wanted.
    \end{itemize}
The lemma follows.
 \hfill $\blacksquare$\\

We can now treat the general case.\\

\noindent \textsc{Proof of Theorem~\ref{minproductorder}.}  We will proceed by induction on $\lwidth(\Delta^{(1)})$. The base case is where $\lwidth(\Delta^{(1)}) \le 0$,
for which the theorem can be verified explicitly using Koelman's classification: see \cite[Thm.~10]{movingout} or \cite[Ch.~4]{Koelman}.

So assume that $\Delta^{(1)}$ is two-dimensional.  
Because $\lwidth(\Delta^{(1)(1)}) < \lwidth(\Delta^{(1)})$ we can apply the induction hypothesis
to find that $\Delta^{(1)}$ can be positioned inside the box
\begin{equation} \label{hypobox}
[0, \lwidth(\Delta^{(1)}) ] \ \times \ [0, \lsize_\square(\Delta^{(1)})  ].
\end{equation}
The foregoing lemma allows us to assume that $\Delta^{(1)}$ is not a standard triangle. But then $\Delta$ must
be contained in the strip
\begin{equation} \label{verticalstrip}
  [-1, \lwidth(\Delta) -1 ] \ \times \ \mathbb{R}
\end{equation}
by Lemma~\ref{lubbesschicholemma}.
  We make a case distinction:
\begin{itemize}
  \item Suppose that the box (\ref{hypobox}) is a square, i.e.\ $\lwidth(\Delta^{(1)}) = \lsize_\square(\Delta^{(1)})$. 
  Then by symmetry
  $\Delta$ must also be contained in the strip
  \[ \mathbb{R} \ \times \ [-1, \lwidth(\Delta) -1 ] . \]
  So it is contained in the intersection
  \[ [-1, \lwidth(\Delta) - 1 ] \ \times \ [-1, \lwidth(\Delta) - 1 ]  \ \simeq \ [0, \lwidth(\Delta) ] \ \times \ [0, \lwidth(\Delta) ]. \]
  Therefore $\lsize_\square(\Delta) = \lwidth(\Delta)$, and the statement follows.
  \item Suppose that the box (\ref{hypobox}) is not a square, i.e.\ $\lwidth(\Delta^{(1)}) < \lsize_\square(\Delta^{(1)})$. 
  We make a further distinction:
  \begin{itemize}
    \item Suppose that an edge of $\Delta^{(1)}$ arises as a slice of $\lsize_\square(\Delta^{(1)}) \cdot \square$.
    Because $\lwidth(\Delta^{(1)}) < \lsize_\square(\Delta^{(1)})$ it necessarily concerns one of the two vertical edges of our box (\ref{hypobox}). By flipping
    horizontally if needed we can assume that it concerns the left edge, which is then a common edge of $\Delta^{(1)}$ with 
    $\lsize_\square(\Delta^{(1)}) \cdot \square$.
    \begin{center}
  \psset{unit=0.3cm}
  \begin{pspicture}(-1,-2.8)(7,7.5)
  \pspolygon[linestyle=dashed](0,0)(4,0)(4,6)(0,6)
  \pspolygon[fillstyle=solid,fillcolor=lightgray](0,0)(2,1)(4,4)(1,6)(0,6)
  \psline[linestyle=dotted](4,0)(6,0)
  \psline[linestyle=dotted](6,0)(6,6)
  \psline[linestyle=dotted](6,6)(4,6)
  \pscurve(0,-0.5)(0.5,-0.8)(2.6,-0.8)(3,-1.2)
  \pscurve(6,-0.5)(5.5,-0.8)(3.4,-0.8)(3,-1.2)
  \pscurve(0,6.5)(0.5,6.8)(1.6,6.8)(2,7.2)
  \pscurve(4,6.5)(3.5,6.8)(2.4,6.8)(2,7.2)
  \rput(3.1,-2.1){\small $\lsize_\square(\Delta^{(1)})$}
  \rput(2.1,8.1){\small $\lwidth(\Delta^{(1)})$}
  \pscurve{<-}(5.8,6.2)(6.5,6.8)(7.2,6.6)
  \rput(10.9,6.6){\small $\lsize_\square(\Delta^{(1)}) \cdot \square$}
  \end{pspicture}
  \end{center}
    It is the unique such edge, so we can apply
    Corollary~\ref{skewingcorollary} to find
    a vertical skewing $\varphi$ such that $\varphi(\Delta) \subset (-1,-1) + \lsize_\square(\Delta) \cdot \square$.
    But the strip (\ref{verticalstrip}) is invariant under vertical skewings. By taking the intersection, we find that
    \[ \varphi(\Delta) \subset [-1, \lwidth(\Delta) - 1 ] \ \times \ [-1, \lsize_\square(\Delta) - 1 ] \ \simeq \ [0, \lwidth(\Delta) ] \ \times \ [0, \lsize_\square(\Delta) ], \]
    as wanted.
    \item Suppose that $\Delta^{(1)}$ does not have an edge arising as a slice of $\lsize_\square(\Delta^{(1)}) \cdot \square$. Because
    our box (\ref{hypobox}) is not a square, it cannot have a diagonal edge either. So from 
    Corollary~\ref{slicingcorollary} we see that 
    \[ \Delta \subset (-1,-1) + (\lsize_\square(\Delta^{(1)}) + 2) \cdot \square =  (-1,-1) + \lsize_\square(\Delta) \cdot \square, \]
    and we conclude as above.
  \end{itemize}
\end{itemize}
The theorem follows. \hfill $\blacksquare$\\
 
We conclude this section by remarking that the above material can be used to design an algorithm for simplifying 
rational parametrizations of toric surfaces, following Schicho~\cite{schicho}, where the focus now lies on the bidegree rather than the total degree.

\section{Proof of the main theorem} \label{section_maintheorem}

After this large chunk of combinatorics, let us return to algebraic geometry. As in the introduction, let $k$ be an algebraically closed field (of arbitrary characteristic), 
let $f \in k[x^{\pm 1}, y^{\pm 1}]$ be an irreducible Laurent polynomial, and assume that $\Delta(f)$ is two-dimensional. 
Let $U_f \subset \mathbb{T}^2$ be the curve defined by $f$. 
Our aim is to prove Theorem~\ref{s2upper}. \\
 
 \noindent \textsc{Proof of Theorem~\ref{s2upper}.} 
 First remark that the inequalities (\ref{tobeprovenupperbounds}) are trivial as soon as $U_f$ is a rational curve, because the right-hand sides 
 are at least $1$ resp.\ $(1,1)$. In particular, by Baker's bound \cite{beelen} we can assume that $\Delta(f)^{(1)}$ is not empty. 
 But then the right-hand sides are at least $3$ and $(2,2)$. For curves of genus one these bounds can be met simultaneously. Indeed, pick a cubic (e.g.\ Weierstrass) model
 and apply a projective transformation ensuring that the curve passes through the two coordinate points at infinity. Then its affine part is defined by a polynomial whose
 Newton polygon is contained in 
   \begin{center}
  \psset{unit=0.5cm}
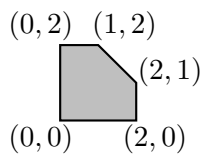
  \begin{pspicture}(-1,-0.7)(3,3)
  \pspolygon[fillstyle=solid,fillcolor=lightgray](0,0)(2,0)(2,1)(1,2)(0,2)
  \rput(-0.5,-0.5){\small $(0,0)$}
  \rput(2.5,-0.5){\small $(2,0)$}
  \rput(2.9,1.3){\small $(2,1)$}
  \rput(1.7,2.5){\small $(1,2)$}
  \rput(-0.5,2.5){\small $(0,2)$}
  \end{pspicture}
  \end{center}
  and therefore both in $3\Sigma$ and $2\square$.
  Thus we can assume that $U_f$ is of genus $g \geq 2$. By Baker's bound this implies that $\sharp(\Delta(f)^{(1)} \cap \mathbb{Z}^2) \geq 2$.
 
 Let us begin with proving the first inequality $s_2(U_f) \leq \lsize_\Sigma(\Delta(f)^{(1)}) + 3$. By the trivial bound
 (\ref{trivialbounds}) it suffices to analyze the exceptional polygons listed in~Theorem~\ref{latticesize_classification}. Since
 $\sharp(\Delta(f)^{(1)} \cap \mathbb{Z}^2) \geq 2$ this leaves us with two cases:
 \begin{itemize}
   \item Assume that $\Delta(f) = \square_{a,b}$ for certain $a,b \geq 2$. Pick a point $(x_0, y_0) \in U_f$. Then the Newton polygon of
 $f(x + x_0, y + y_0)$ is contained in 
 \[ \conv \{ (1,0), (a,0), (a, b), (0, b), (0,1) \}. \]
 But then $x^ay^b f(x^{-1} + x_0, y^{-1} + y_0)$ is a polynomial of degree at most $a + b - 1$.
 So $s_2(U_f) \leq \lsize_\Sigma(\Delta(f)) - 1 = \lsize_\Sigma(\Delta(f)^{(1)}) + 3$.
\item Assume that there exist parallel edges
 $\tau \subset \Delta(f)$ and $\tau' \subset \Delta(f)^{(1)}$ whose
supporting lines are at integral distance $1$ of each other, 
of respective lengths $r$ and $s$, such that $r \geq s + 4$. From Lemma~\ref{latticesize_auxiliary1} and its proof we see
that $s = \lsize_\Sigma(\Delta(f)^{(1)})$ and that we can 
assume that 
\[ \tau = \conv \{ (0,0), (r,0) \} \qquad \text{and} \qquad \tau' = \conv \{ (1,1), (s + 1, 1) \}. \] This configuration implies
that $\Delta(f)$ is contained in the half-planes $X \geq 0$, $Y \geq 0$ and $X + (r - s - 2)Y \leq r$.
In other words, 
\[ f = \sum_{j = 0}^{\lfloor r / (r - s - 2) \rfloor} g_j(x) y^j \]
for polynomials $g_j \in k[x]$ satisfying $\deg g_j \leq r - (r - s - 2)j$ and $\deg g_0 = r$. Now factor
$g_0(x) = g_0'(x) h_0'(x)$ with $\deg g_0' = s + 3$ and $\deg h_0' = r - s - 3$, substitute $y \leftarrow y h_0'(x)$, and kill a factor $h_0'(x)$ to obtain
\[ g_0'(x) + \sum_{j = 1}^{\lfloor r / (r - s - 2) \rfloor}  g_j(x)h_0'(x)^{j-1} y^j. \]
One verifies that each term has degree at most $s + 3$, which proves that $s_2(U_f) \leq s + 3 = \lsize_\Sigma(\Delta(f)^{(1)}) + 3$.
 \end{itemize}
As for the case where $\Delta(f) \simeq d\Upsilon$
for some $d \geq 2$, note that by Theorem~\ref{latticesize_classification} we have $\lsize_\Sigma(d\Upsilon) = 3d$, so the bound
we need to prove is sharper. Consider the embedding
\[ \psi:  \mathbb{T}^2 \hookrightarrow \mathbb{P}^3 = \text{Proj} \, [X_0, X_1, X_2, X_3] : (x,y) \mapsto (x^{-1}y^{-1} : x : y : 1). \]
It embeds $U_f$ in a projective curve $C_f$ which arises as the intersection of the cubic $X_0X_1X_2 - X_3^3$ and an irreducible hypersurface of degree $d$, whose concrete 
equation depends on $f$. In particular it is a curve of degree $3d$.
By \cite[IV.Prop.\ 3.8 and IV.Thm.\ 3.9]{hartshorne} we can find a point on $C_f$, 
the general secant line through which is not a multisecant.
Projecting from such a point yields
a birational equivalence between $C_f$ and 
a plane curve of degree $3d-1$, as wanted. 

 Next we address the inequality 
 \[ s_{1,1}(U_f) \leq s_{1,1}(\Delta(f)^{(1)}) + (2,2) = \left( \lwidth(\Delta(f)^{(1)}) + 2, \lsize_\square(\Delta(f)^{(1)}) + 2 \right). \] 
 We make a case distinction. 
 \begin{itemize}
 \item Assume
 that $\Delta(f)^{(1)} \simeq (d-3) \Sigma$ for some $d \geq 4$, so that 
 \[ \lsize_\square(\Delta(f)^{(1)}) = \lsize_\Sigma(\Delta(f)^{(1)}) = d - 3. \] By the foregoing $U_f$ has a plane model of
 degree $d$. Using a projective transformation we can ensure that this model passes through the coordinate points at infinity. As
 in the genus one case we end up
 with a model of bidegree $(d-1,d-1)$, as wanted. 
 \item Suppose that $\Delta(f)^{(1)}$ is not a standard triangle. By Lemma~\ref{lubbesschicholemma} we have $ \lwidth(\Delta(f)) = \lwidth(\Delta(f)^{(1)}) + 2$.
 If we are not among the exceptions listed in Theorem~\ref{classificationsquare} then also 
 $\lsize_\square(\Delta(f)) = \lsize_\square(\Delta(f)^{(1)}) + 2$, and the statement follows from the bound (\ref{trivialbounds}). 
 
 Because
 $\sharp(\Delta(f)^{(1)} \cap \mathbb{Z}^2) \geq 2$ only the last exception is a concern.
 Assume that there exist parallel edges
 $\tau \subset \Delta(f)$ and $\tau' \subset \Delta(f)^{(1)}$ whose
supporting lines are at integral distance $1$ of each other, 
of respective lengths $r$ and $s$, such that $r \geq s + 3$. 
By Lemma~\ref{latticesizesquare_auxiliary1} we know that $s = \lsize_\square(\Delta(f)^{(1)})$. Thus our aim
is to apply a birational change of variables to $f$ so that the result has bidegree $(\lwidth(\Delta(f)), s+2)$.

Again, as in the proof of Lemma~\ref{latticesize_auxiliary1} we can assume that 
 \[ \tau = \conv \{ (0,0), (r,0) \} \qquad \text{and} \qquad \tau' = \conv \{ (1,1), (s + 1, 1) \}, \]
so that $\Delta(f)$ is contained in the half-planes $X \geq 0$, $Y \geq 0$ and $X + (r - s - 2)Y \leq r$.
This implies that $\Delta(f)^{(1)}$ is contained in $(1,1) + s \Sigma$, and because we excluded standard triangles the top
vertex of the latter cannot occur, from which
one sees that $\lwidth(\Delta(f)^{(1)}) < s$. 

If we now use Theorem~\ref{minproductorder} to position $\Delta(f)^{(1)}$ inside a box
\[ [ 1, s + 1 ] \times [ 1, \lwidth(\Delta(f)^{(1)}) + 1], \]
then $\tau'$ necessarily arises as a horizontal line segment; we can assume it to be
the bottom segment
$\conv \{ (1,1), (s + 1, 1) \}$. By Lemma~\ref{lubbesschicholemma} our Newton polygon $\Delta(f)$ is then contained
in the strip 
\begin{equation} \label{horizontalstripje}
 \mathbb{R} \ \times \ [ 0, \lwidth(\Delta(f))].
\end{equation} 
Now once again as in the proof of Lemma~\ref{latticesize_auxiliary1}
we can apply a horizontal skewing to position $\tau$ at $\conv \{ (0,0), (r,0) \}$. 
We again obtain that $\Delta(f)$ is contained in the half-planes $X \geq 0$, $Y \geq 0$ and $X + (r - s - 2)Y \leq r$, while it is also
kept in the strip (\ref{horizontalstripje}).
In other words, 
\[ f = \sum_{j = 0}^{\lwidth(\Delta(f))} g_j(x) y^j \]
for polynomials $g_j \in k[x]$ satisfying $\deg g_j \leq r - (r - s - 2)j$, $\deg g_0 = r$ and $g_{\lwidth(\Delta(f))} \neq 0$. 
Now factor
$g_0(x) = g_0'(x) h_0'(x)$ with $\deg g_0' = s + 2$ and $\deg h_0' = r - s - 2$, substitute $y \leftarrow y h_0'(x)$, and kill a factor $h_0'(x)$ to obtain
a polynomial
\[ g_0'(x) + \sum_{j = 1}^{\lwidth(\Delta(f))}  g_j(x)h_0'(x)^{j-1} y^j \]
of degree $s+2$ in $x$ and degree $\lwidth(\Delta(f)))$ in $y$, as wanted. 
   \end{itemize}
It remains to show that $s_{1,1}(U_f) \leq (3,4)$ when $\Delta(f) = 2\Upsilon$. 
By Baker's bound $U_f$ is a curve of genus at most $4$. If $U_f$ is hyperelliptic then
the bound follows trivially (because of the lexicographic order). If $U_f$ is non-hyperelliptic
of genus $3$ then $U_f$ is birationally equivalent to a non-singular quartic in $\mathbb{P}^2$,
and one can construct a model of bidegree $(3,3)$ by forcing it through the coordinate points at infinity.
Finally if $U_f$ is non-hyperelliptic of genus $4$ then it is birationally equivalent to a singular quintic in $\mathbb{P}^2$ by 
 \cite[IV.Ex.\,5.4]{hartshorne}. Using a projective transformation we can assume that the curve passes through the coordinate
 points at infinity, one of these being a singularity. Dehomogenizing yields an affine model of bidegree $(3,4)$ as wanted.
\hfill $\blacksquare$

\section{Cases where the bounds are sharp} \label{section_conjecture}

In this section we again restrict to $\text{char} \, k = 0$, because of some references on which we will rely.
One of these references is a subsequent, more elaborate paper \cite{linearpencils} of ours, in which we study  
linear pencils that are encoded in the Newton polygon. At some point in that paper, the lattice size with respect to $\Sigma$
pops up as a convenient notion \cite[Thm.~7.2]{linearpencils}. This is how we came up with the first inequality from Theorem~\ref{s2upper},
which meant the start of this project.

We will make extensive reference to \cite{linearpencils}, even though it concerns a successive paper. But we stress that no circular
reasoning is being made: no statements in \cite{linearpencils}
make use of any of the results of this section. Moreover, some of the results of \cite{linearpencils} that we need
appear (in more disguised terms) in an earlier article by Kawaguchi \cite{kawaguchi}. Finally, we emphasize
that the primary aim of this section is to convince the reader that 
the bounds from Theorem~\ref{s2upper} often give the correct values of $s_2(U_f)$ and $s_{1,1}(U_f)$, and to give
some evidence in favor of Conjecture~\ref{sufficientlygenericconjecture}; we will not push the limits of our exposition.

Let us specify what we mean by $f$ being sufficiently generic with respect to its Newton polygon $\Delta(f)$.
To each two-dimensional lattice polygon $\Delta$ there is a standard way of associating
a toric surface $X(\Delta)$ over $k$ (along with an embedding in $\mathbb{P}^{\# (\Delta \cap \mathbb{Z}^2) - 1}$). This is a completion of the torus $\mathbb{T}^2$, so it is natural
to consider the closure of $U_f$ inside it. 
It turns out that for almost all Laurent polynomials $f$ the closure $C_f$ of $U_f$ inside
$X(\Delta(f))$ is non-singular. More precisely, if one fixes a two-dimensional lattice polygon $\Delta$, then the locus of the Laurent polynomials
$f$ for which $\Delta(f) = \Delta$ and $C_f$ is non-singular is dense in the according $\#(\Delta \cap \mathbb{Z}^2)$-dimensional
coefficient space. We refer to \cite[\S2]{castryckvoight} and \cite[\S4]{linearpencils} for more background.

We now rephrase Conjecture~\ref{sufficientlygenericconjecture} as follows.

\begin{conjecture} \label{theguess}
 If $C_f$ is a non-singular curve and $\Delta(f) \not \simeq 2\Upsilon$ then 
 \[ s_2(U_f) = \lsize_\Sigma(\Delta(f)^{(1)}) + 3 \quad \text{and} \quad s_{1,1}(U_f) = s_{1,1}(\Delta(f)^{(1)}) + (2,2), \] unless
 $\Delta(f) \simeq d\Upsilon$ for some $d \geq 3$, in which case $s_2(U_f) = 3d - 1$.
\end{conjecture}

\noindent This would extend the list of geometric invariants that are known to be encoded in the Newton polygon. We mention some
of its current entries: if $C_f$ is a non-singular curve then
\begin{itemize}
  \item[(i)] its (geometric) genus $g$ equals $\sharp (\Delta(f)^{(1)} \cap \mathbb{Z}^2)$; this is due to Khovanskii \cite{khovanskii};
  \item[(ii)] its gonality $c$ equals $\lwidth(\Delta(f)^{(1)}) + 2$, unless $\Delta(f) \simeq 2 \Upsilon$ in which case the gonality equals $3$; this
  is \cite[Cor.\,6.2]{linearpencils}, whose proof strongly builds on previous work of Kawaguchi \cite{kawaguchi};
  \item[(iii)] it is isomorphic to a non-singular plane curve if and only if $\Delta(f)^{(1)} = \emptyset$ or $\Delta(f)^{(1)} \simeq (d-3)\Sigma$ for some $d \geq 3$; this is \cite[Cor.\,8.2]{linearpencils}.
\end{itemize}
For an extension of this list we refer to \cite{schreyersinvariants,linearpencils,intrinsicness}. Note the similarity between statement (ii) and Conjecture~\ref{theguess}. 

A moral reason for the fact that many invariants are encoded in the Newton polygon
is that $C_f$ canonically embeds inside $X(\Delta^{(1)}) \subset \mathbb{P}^{g-1}$, and that the defining equations of the latter are
so special (quadrics of very low rank) that they can often be recovered from the canonical ideal of $C_f$ itself. We refer to the introduction of \cite{schreyersinvariants}
for an extended discussion. Up to equivalence, the polygon $\Upsilon$ is the unique two-dimensional polygon of the form $\Delta^{(1)}$ for
which the ideal of $X(\Delta^{(1)})$ is \emph{not} generated by quadrics. This explains the special role of $2\Upsilon$, which is the only polygon having $\Upsilon$ as its interior. 
If $\Delta(f)^{(1)} \simeq 2 \Upsilon$ then 
from (i) we find that $C_f$ is a genus four curve, for which 
\begin{itemize}
\item $s_2(U_f) = 5$, so the formula $s_2(U_f) = 3d - 1$ is actually correct here: the existence of a degree $5$ model follows from
Theorem~\ref{s2upper}, while degree $4$ or less would contradict that the genus is $5$; 
\item $s_{1,1}(U_f) = (3,4)$ or $s_{1,1}(U_f) = (3,3)$,
depending on whether the unique quadric in which $C_f$ canonically embeds is singular or not: in this 
case by \cite[\S6]{castryckvoight} there exists an $f' \in k[x^{\pm 1}, y^{\pm 1}]$ with
\[ \Delta(f') = \conv \{ (0,0), (6,0), (0,3) \} \quad \text{resp.} \quad \Delta(f') = [0,3] \times [0,3], \]
such that $C_{f'}$ is non-singular and birationally equivalent to
$U_f$; the formulas then follow from Theorem~\ref{theguessproofsquare} below.
\end{itemize}
Alternatively, these formulas can be proved along the lines of \cite[IV.Ex.\,5.4]{hartshorne}.

Note that by Lemma~\ref{spectrum} and (\ref{latticesizechain}) we have 
\[ -1 \leq \lwidth(\Delta(f)^{(1)}) \leq \lsize_\Sigma(\Delta(f)^{(1)}) \quad \text{and} \quad
 -1 \leq \lwidth(\Delta(f)^{(1)}) \leq \lsize_\square(\Delta(f)^{(1)}). \] 
We can prove Conjecture~\ref{theguess}
near both ends of these ranges.

\begin{theorem} \label{theguessproof}
 If $C_f$ is non-singular and 
 \[ \lwidth(\Delta(f)^{(1)}) \leq 1 \quad \text{or} \quad \lwidth(\Delta(f)^{(1)}) \geq \lsize_\Sigma(\Delta(f)^{(1)}) - 2, \] then
 $s_2(U_f) =  \lsize_\Sigma(\Delta(f)^{(1)}) + 3$, unless $\Delta(f) \simeq 2\Upsilon, 3\Upsilon$, in which case
 $s_2(U_f) = 5$ and $s_2(U_f) = 8$, respectively.
\end{theorem}

\noindent \textsc{Proof.} 
By the above discussion we can assume that $\Delta(f) \not \simeq 2\Upsilon$.

At the lower end, we can argue as follows:
\begin{itemize}
 \item If $\lwidth(\Delta(f)^{(1)}) = -1$, or in other words if $\Delta(f)^{(1)} = \emptyset$, then $U_f$ is rational because of (i), and there is nothing to prove.
 \item
If $\lwidth(\Delta(f)^{(1)}) = 0$, then $\Delta(f)^{(1)}$ is a line segment, say of integral length
$g-1$. By (i) and (ii) we find that $U_f$ is hyperelliptic of genus $g$. So
$s_2(U_f) = g + 2$, which indeed equals $\lsize_\Sigma(\Delta(f)^{(1)}) + 3$.
 \item If $\lwidth(\Delta(f)^{(1)}) = 1$ then $\Delta(f)^{(1)}$ is equivalent to a Lawrence prism
  \begin{center}
  \psset{unit=0.25cm}
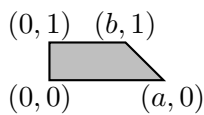
  \begin{pspicture}(-1,-1.2)(7,3.8)
  \pspolygon[fillstyle=solid,fillcolor=lightgray](0,0)(6,0)(4,2)(0,2)
  \rput(-0.5,-1){\small $(0,0)$}
  \rput(6.5,-1){\small $(a,0)$}  
  \rput(-0.5,2.9){\small $(0,1)$}  
  \rput(4,2.9){\small $(b,1)$}  
  \end{pspicture}
  \end{center}
with $1 \leq a \geq b \geq 0$. In this case $C_f$ is a trigonal curve with scrollar invariants $a,b$ 
by \cite[Thm.~9.1]{linearpencils}. From \cite[Lem.~2.1]{KeemMartens}, which is expressed
in terms of the Maroni invariant $b = \min(a,b)$, we conclude
that $s_2(U_f) = g + 1 - b = a + 3$ if $a > b$, and $s_2(U_f) = g + 2 - b = a + 4$ if $a = b$. 
By Theorem~\ref{latticesize_classification}, in both cases this exactly matches with $\lsize_\Sigma(\Delta(f)^{(1)}) + 3$.
\end{itemize}
At the other end, we make the following reasonings. 
\begin{itemize}
\item Assume that $\lwidth(\Delta(f)^{(1)}) = \lsize_\Sigma(\Delta(f)^{(1)})$. Then by
Lemma~\ref{spectrum} we have that $\Delta(f)^{(1)} \cong (d-3)\Sigma$
for some integer $d \geq 3$. But then by (iii) our curve $C_f$ is isomorphic to a non-singular plane curve of degree $d$, and therefore $s_2(U_f) = d = \lsize_\Sigma(\Delta(f)^{(1)}) + 3$.
\item
If $\lwidth(\Delta(f)^{(1)}) = \lsize_\Sigma(\Delta(f)^{(1)}) - 1$ then by (ii) the gonality of $U_f$ equals
\[ c = \lsize_\Sigma(\Delta(f)^{(1)}) + 1, \]
unless $\Delta(f) \simeq 2 \Upsilon$, but this case was excluded. On the other hand, again by (iii) every plane model is necessarily singular. 
This means that $s_2(U_f) \geq c+2$, because otherwise
projection from a singular point on the plane model would give a map to $\mathbb{P}^1$ of degree strictly less than $c$. We conclude that
$s_2(U_f) \geq \lsize_\Sigma(\Delta(f)^{(1)}) + 3$, and by Theorem~\ref{s2upper} equality holds.
\item If $\lwidth(\Delta(f)^{(1)}) = \lsize(\Delta(f)^{(1)}) - 2$ then by (ii) the gonality of $U_f$ equals
\[ c = \lsize_\Sigma(\Delta(f)^{(1)}). \]
By (iii) every plane model is singular and we find as above that $s_2(U_f) \geq c+2$. In case $\Delta(f) \simeq 3 \Upsilon$
this matches with the upper bound from Theorem~\ref{s2upper}, and we are done. We would like to show
that $s_2(U_f) \geq c+3$ in the other cases. So suppose
that $\Delta(f) \not \simeq 3 \Upsilon$ and assume by contradiction that $s_2(U_f) = c + 2$. In this case we see that
the curve carries infinitely many base-point free $g^1_{c+1}$'s, obtained by projection from
the non-singular points of a plane degree $c+2$ model. By
\cite[Thm.~7.2]{linearpencils} this is possible only if 
 $\Delta(f)^{(1)} \cong \conv \{ (0,0), (1,0), (3,1), (2,2), (1,2) \}$, so that $\Delta(f)$ is of the form
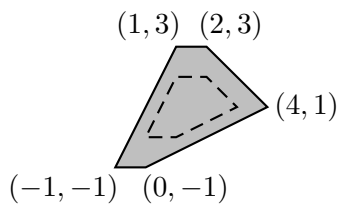
\begin{center}
\psset{unit=0.4cm}
\begin{pspicture}(-1,-1.9)(4,4)
\pspolygon[fillstyle=solid,fillcolor=lightgray](-1,-1)(0,-1)(4,1)(2,3)(1,3)
\pspolygon[linestyle=dashed](0,0)(1,0)(3,1)(2,2)(1,2)
\rput(-2.7,-1.7){\small $(-1,-1)$}
\rput(1.3,-1.7){\small $(0,-1)$}
\rput(5.3,1){\small $(4,1)$}
\rput(2.8,3.7){\small $(2,3)$}
\rput(0.1,3.7){\small $(1,3)$}
\end{pspicture}
\end{center}
(the dashed polygon indicates $\Delta(f)^{(1)}$).
By (i) and (ii) our curve $C_f$ has gonality $c=4$ and geometric genus $g=7$. By \cite[Cor.\ 6.3]{linearpencils} and \cite[Thm.\ 9.1]{linearpencils} the gonality pencil
is unique and its scrollar invariants are $1$, $1$, and $2$. Now take a curve in $\mathbb{P}^2$ of degree $d = c+2 = 6$ that is birationally
equivalent to $C_f$.
Because the gonality is $4$ and the gonality pencil is unique, 
the curve must have a unique singular point $P$, of multiplicity $2$.
The point $P$ cannot be an ordinary node or a cusp, otherwise the genus would be $(d-1)(d-2)/2 - 1 = 9$.
Thus there is a unique tangent line at $P$ which intersects
the curve at $P$ with multiplicity at least $4$.
Using a transformation of $\mathbb{P}^2$ we can assume that $P = (0:1:0)$ and that this line is at infinity. Dehomogenizing
the defining equation then results in a polynomial $g(x,y)$ that is supported on the following polygon:
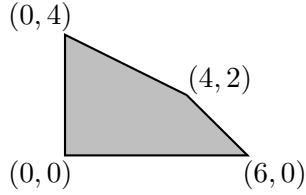
\begin{center}
\psset{unit=0.4cm}
\begin{pspicture}(0,-1)(6,5)
\pspolygon[fillstyle=solid,fillcolor=lightgray](0,0)(6,0)(4,2)(0,4)
\rput(-0.8,-0.65){\small $(0,0)$}
\rput(6.9,-0.65){\small $(6,0)$}
\rput(5.1,2.5){\small $(4,2)$}
\rput(-0.8,4.6){\small $(0,4)$}
\end{pspicture}
\end{center}
The coefficient at $y^4$ is non-zero because the gonality is $4$. In particular
our $g^1_4$ is given by the projection $(x,y) \mapsto x$. 
Also note that at least one of the coefficients
at $x^4y^2$, $x^5y$, $x^6$ is non-zero, because the degree is $6$.
Now let $D_x \in g^1_4$ be the zero divisor of
$x^{-1}$, and similarly let $D_y$ be the zero divisor of $y^{-1}$. The steepness of the above polygon ensures that $D_{y} \leq 2D_{x}$. In particular
$H^0(2D_{x}) \supset \{ 1, x^{-1}, x^{-2}, y^{-1} \}$ is at least $4$-dimensional. This shows that $0$ must be among the scrollar invariants of our $g^1_4$: a contradiction.
We conclude that $s_2(U_f) \geq c + 3$: a contradiction. \hfill $\blacksquare$\\
\end{itemize}

\begin{theorem} \label{theguessproofsquare}
 If $C_f$ is non-singular, $\Delta(f) \not \simeq 2 \Upsilon$ and
 \[ \lwidth(\Delta(f)^{(1)}) \leq 1 \quad \text{or} \quad \lwidth(\Delta(f)^{(1)}) \geq \lsize_\square(\Delta(f)^{(1)}) - 1, \] then
 $s_{1,1}(U_f) =  s_{1,1}(\Delta(f)^{(1)}) + (2,2)$.
\end{theorem}

\noindent \textsc{Proof.} The case of rational and (hyper)elliptic curves follows from Lemma~\ref{marclemma}. In the trigonal case:
\begin{itemize}
\item If $g = 3$ then $\Delta(f)^{(1)} \simeq \Sigma$ and the upper bound from Theorem~\ref{s2upper} reads $(3,3)$, which
is clearly sharp in the case of a trigonal curve.
\item If $g = 4$ then either $\Delta(f)^{(1)} \simeq \conv \{ (0,0), (2,0), (0,1) \}$ or $\Delta(f)^{(1)} \simeq \square_{1,1}$; indeed $\Delta(f)^{(1)} \simeq \Upsilon$ 
was excluded in the statement of the theorem. In the latter case the upper bound from Theorem~\ref{s2upper} reads $(3,3)$, which is clearly optimal in the case of a trigonal curve.
In the former case the upper bound reads $(3,4)$, which is also optimal because by \cite[Cor.\,6.3]{linearpencils} the gonality pencil is unique, while the
existence of a model of bidegree $(3,3)$ would contradict that.
\item If $g \geq 5$ then the $g^1_3$ is always unique. From
\cite[Prop.\,1]{martensschreyer} (see also \cite[Ex.\,1.2.7]{CoppensKeemMartens}) one sees that there exists a base-point free $g^1_d$ that is independent from
this $g^1_3$ if and only if $g - d$ does not exceed the Maroni invariant. Using the same notation as in the foregoing proof, this condition reads $g -d \leq a$, which is equivalent
with $d \geq b + 2$. Thus $s_{1,1}(U_f) = (3,b+2) = s_{1,1}(\Delta(f)) + (2,2)$, as wanted.
\end{itemize}
At the other end, we make the following reasonings.
\begin{itemize}
  \item If $\lwidth(\Delta(f)^{(1)}) = \lsize_\square(\Delta(f)^{(1)})$ then by (ii) the gonality of $U_f$ equals
  \[ c = \lwidth(\Delta(f)^{(1)}) + 2 = \lsize_\square(\Delta(f)^{(1)}) + 2 \]
  unless $\Delta(f) \cong 2\Upsilon$, but this case was excluded. 
  So $s_{1,1}(\Delta(f)^{(1)}) + (2,2) = (c,c)$ is clearly a lower bound for $s_{1,1}(U_f)$,
  and by Theorem~\ref{s2upper} equality holds.
  \item If $\lwidth(\Delta(f)^{(1)}) = \lsize_\square(\Delta(f)^{(1)}) - 1$
  then  \[ c = \lwidth(\Delta(f)^{(1)}) + 2 = \lsize_\square(\Delta(f)^{(1)}) + 1, \]
  and we similarly find that $s_{1,1}(\Delta(f)^{(1)}) + (2,2) = (c,c + 1)$. So it
  is sufficient to show that $(c,c) \notin S_{1,1}(U_f)$, i.e.\ our curve does not carry two independent
  gonality pencils. But by \cite[Thm.\,6.1]{linearpencils} every gonality pencil is combinatorial, i.e.\ it corresponds to
  projecting along some lattice width direction of $\Delta(f)$. In particular $\lwidth(\Delta(f)) = c$, and if
  $C_f$ would admit two gonality pencils then $\Delta(f)$ would admit two $\mathbb{R}$-linearly independent
  lattice width directions. By \cite[Lem.\,5.2(v)]{linearpencils} this would mean that $\lsize_\square(\Delta(f)) = c$. But then 
  \[ c - 1 = \lsize_\square(\Delta(f)^{(1)}) \leq c-2,\] a
  contradiction.\hfill $\blacksquare$
\end{itemize}

There is room for improvement in 
Theorems~\ref{theguessproof} and~\ref{theguessproofsquare}, in order to cover larger ranges of $\text{lw}(\Delta(f)^{(1)})$.
At the lower end this seems difficult however. Whereas it is well-understood which base-point free pencils occur in the hyperelliptic and trigonal cases~\cite{martensschreyer},
for curves of higher gonality not much seems known, although Coppens and Martens proved some potentially useful facts in the tetragonal case \cite{CoppensMartens}.
At the upper end more seems possible: one can try to extend the results of \cite[\S7]{linearpencils} in order to describe the $g^1_{c + n}$'s on smooth curves
in toric surfaces, for $n=2,3, \dots$\ It is expected that these are always combinatorially determined, except for a finite (but increasing) number of polygons. This 
would help in pushing the above arguments. The
finite number of exceptions can then hopefully be treated using an ad hoc idea, such as the one used in the proof of Theorem~\ref{theguessproof}. We expect this to
become increasingly difficult and case-distinctive, however. Alternatively, it might be possible to obtain some results by using specialization of linear systems from curves to graphs \cite{baker,curvegonalities} to reduce Conjecture~\ref{sufficientlygenericconjecture} to a purely combinatorial statement.

\small

\noindent \textsc{Vakgroep Wiskunde}\\ 
\noindent \textsc{Universiteit Gent}\\
\noindent \textsc{Krijgslaan 281, 9000 Gent, Belgium}\\
\noindent \emph{E-mail address:} \verb"wouter.castryck@gmail.com"\\

\noindent \textsc{Department of Mathematics and Applied Mathematics}\\ 
\noindent \textsc{University of Cape Town}\\
\noindent \textsc{Private Bag X1, Rondebosch 7701, South Africa}\\
\noindent \emph{E-mail address:} \verb"filip.cools@uct.ac.za"\\

\end{document}